# OPTIMAL BOUNDS FOR AGGREGATION OF AFFINE ESTIMATORS

BY PIERRE C. BELLEC[1]

*ENSAE and Rutgers University*

We study the problem of aggregation of estimators when the estimators are not independent of the data used for aggregation and no sample splitting is allowed. If the estimators are deterministic vectors, it is well known that the minimax rate of aggregation is of order $\log(M)$, where $M$ is the number of estimators to aggregate. It is proved that for affine estimators, the minimax rate of aggregation is unchanged: it is possible to handle the linear dependence between the affine estimators and the data used for aggregation at no extra cost. The minimax rate is not impacted either by the variance of the affine estimators, or any other measure of their statistical complexity. The minimax rate is attained with a penalized procedure over the convex hull of the estimators, for a penalty that is inspired from the $Q$-aggregation procedure. The results follow from the interplay between the penalty, strong convexity and concentration.

**1. Introduction.** We study the problem of recovering an unknown vector $\mathbf{f} = (f_1, \dots, f_n)^T \in \mathbf{R}^n$ from noisy observations:

$$Y_i = f_i + \xi_i, \qquad i = 1, \dots, n, \tag{1.1}$$

where the noise random variables $\xi_1, \dots, \xi_n$ are i.i.d. $\mathcal{N}(0, \sigma^2)$ or i.i.d. sub-Gaussian random variables. We measure the quality of estimation of the unknown vector $\mathbf{f}$ with the squared Euclidean norm in $\mathbf{R}^n$:

$$\|\mathbf{f} - \hat{\boldsymbol{\mu}}\|_2^2,$$

for any estimator $\hat{\boldsymbol{\mu}}$ of $\mathbf{f}$. When the noise random variables are normal, (1.1) is the Gaussian sequence model, which has been extensively studied; see, for example, [22] and the references therein. Several estimators have been proposed to recover the unknown vector $\mathbf{f}$ from the observations: the ordinary least squares, the ridge estimator, the Stein estimator and the procedures based on shrinkage, to name a few. Several of these estimators depend on a parameter that must be chosen carefully to obtain satisfying error bounds. These available estimators have different

Received September 2015; revised December 2016.

[1]Supported by GENES and by the grant Investissements d'Avenir (ANR-11-IDEX-0003/Labex Ecodec/ANR-11-LABX-0047).

[MSC2010 subject classifications.](#) Primary 62G05; secondary 62J07.

strengths and weaknesses in different scenarios, so it is important to be able to mimic the best among a given family of estimators, without any assumption on the unknown **f**. The problem of mimicking the best estimator in a given finite set is the problem of model-selection type aggregation, which was introduced in [30, 37]. More precisely, let $\hat{\boldsymbol{\mu}}_1, \ldots, \hat{\boldsymbol{\mu}}_M$ be $M$ estimators of **f** based on the data $\mathbf{y} = (Y_1, \ldots, Y_n)^T$. The goal is to construct with the same data $\mathbf{y} = (Y_1, \ldots, Y_n)^T$ a new estimator $\hat{\boldsymbol{\mu}}$ called the aggregate, which satisfies with probability greater than $1 - \delta$ the sharp oracle inequality[2]

$$(1.2) \qquad \|\hat{\boldsymbol{\mu}} - \mathbf{f}\|_2^2 \le \min_{j=1,\ldots,M} \|\hat{\boldsymbol{\mu}}_j - \mathbf{f}\|_2^2 + \text{PRICE}_M(\delta),$$

where $\text{PRICE}_M(\cdot)$ is a function of $\delta$ that should be small. The term $\text{PRICE}_M(\cdot)$ will be referred to as the price to pay for aggregating the estimators $\hat{\boldsymbol{\mu}}_1, \ldots, \hat{\boldsymbol{\mu}}_M$. If the estimators $\hat{\boldsymbol{\mu}}_1, \ldots, \hat{\boldsymbol{\mu}}_M$ are deterministic vectors, the price to pay for aggregating these estimators is of order $\sigma^2 \log(M/\delta)$ and (1.2) is satisfied for an estimator $\hat{\boldsymbol{\mu}}$ based on $Q$-aggregation [11]. Considering deterministic estimators is of interest if two independent samples are available, so that $\hat{\boldsymbol{\mu}}_1, \ldots, \hat{\boldsymbol{\mu}}_M$ are based on the first sample while aggregation is performed using the second sample. Then the first sample can be considered as frozen at the aggregation step (for more details see [39]). If the estimators are random (dependent on the data **y** used for aggregation), two natural questions arise:

1. Does the price to pay for aggregation increase because of the dependence between $\hat{\boldsymbol{\mu}}_1, \ldots, \hat{\boldsymbol{\mu}}_M$ and the data **y**, or is it still of order $\sigma^2 \log(M/\delta)$? Is there an extra price to pay to handle the dependence?

2. A natural quantity that captures the statistical complexity of a given estimator $\hat{\boldsymbol{\mu}}_j$ is the variance defined by $\mathbb{E}\|\hat{\boldsymbol{\mu}}_j - \mathbb{E}\hat{\boldsymbol{\mu}}_j\|_2^2$. When the estimators are deterministic, their variances are all zero. Now that the estimators are random, does the price to pay for aggregation depend on the statistical complexities of the estimators $\hat{\boldsymbol{\mu}}_1, \ldots, \hat{\boldsymbol{\mu}}_M$, for example, through their variances? Is it harder to aggregate estimators with large statistical complexities?

The goal of this paper is to answer these questions for affine estimators.

Among the procedures available to estimate **f**, several are linear in the observations $Y_1, \ldots, Y_n$. It is the case for the least squares and the ridge estimators, whereas the shrinkage estimators and the Stein estimator are nonlinear functions of the observations. Examples of estimators that are linear or affine in the observations is given in [12], Section 1.2, [1] and references therein. An affine estimator is of the form $\hat{\boldsymbol{\mu}}_j = A_j \mathbf{y} + \boldsymbol{b}_j$ for a deterministic matrix $A_j$ of size $n \times n$ and a deterministic vector $\boldsymbol{b}_j \in \mathbf{R}^n$. The linearity of the estimators $\hat{\boldsymbol{\mu}}_1, \ldots, \hat{\boldsymbol{\mu}}_M$ makes it possible to explicitly treat the dependence between the estimators $\hat{\boldsymbol{\mu}}_1, \ldots, \hat{\boldsymbol{\mu}}_M$ and the data $\mathbf{y} = (Y_1, \ldots, Y_n)^T$ used to aggregate them. Donoho et al. [14] proved that for orthosymmetric quadratically convex sets (which include all ellipsoids and hyperrectangles), the minimax risk among all linear estimators is within 25% of the minimax risk among all estimators.

---

[2] By sharp, we mean that the constant in front of the term $\min_{j=1,\ldots,M} \|\hat{\boldsymbol{\mu}}_j - \mathbf{f}\|_2^2$ is 1.

The papers [10, 12, 27] derived different procedures that satisfy sharp oracle inequalities for the problem of aggregation of affine estimators when the noise random variables are Gaussian. Leung and Barron [27], Dalalyan and Salmon [12] proposed an estimator $\hat{\boldsymbol{\mu}}^{\mathrm{EW}}$ based on exponential weights, for which the following sharp oracle inequality holds in expectation:

$$\mathbb{E}\|\mathbf{f} - \hat{\boldsymbol{\mu}}^{\mathrm{EW}}\|_2^2 \le \min_{j=1,\ldots,M} \mathbb{E}\|\hat{\boldsymbol{\mu}}_j - \mathbf{f}\|_2^2 + 8\sigma^2 \log M,$$

under the assumption that all $A_j$ are orthoprojectors [orthogonal projection matrices; cf. (1.4)], or under a strong commutativity assumption on the matrices $A_j$. The constant 8 can be reduced to 4 if all $A_j$ are orthoprojectors. If the matrices $A_j$ are not symmetric, [12] achieved a similar oracle inequality by symmetrizing the affine estimators before the aggregation step, which suggests that the symmetry assumption can be relaxed. Although the estimator $\hat{\boldsymbol{\mu}}^{\mathrm{EW}}$ achieves this inequality in expectation, it was shown in [2, 11] that it cannot achieve a similar result in deviation, with an unavoidable error term of order $\sqrt{n}$. In Dai et al. [10], a sharp oracle inequality in deviation is derived for an estimator $\hat{\boldsymbol{\mu}}^Q$ based on $Q$-aggregation [11, 31]. Namely, [10] proves that if the matrices $A_1, \ldots, A_M$ are symmetric and positive semidefinite, the estimator $\hat{\boldsymbol{\mu}}^Q$ satisfies with probability greater than $1 - \delta$:

$$(1.3) \qquad \|\mathbf{f} - \hat{\boldsymbol{\mu}}^Q\|_2^2 \le \min_{j=1,\ldots,M} (\|\hat{\boldsymbol{\mu}}_j - \mathbf{f}\|_2^2 + 4\sigma^2 \operatorname{Tr}(A_j)) + C\sigma^2 \log(M/\delta),$$

where the constant $C$ is proportional to the largest operator norm of the matrices $A_1, \ldots, A_M$. The term $4\sigma^2 \operatorname{Tr}(A_j)$ is intimately linked to the statistical complexity of the estimator $\hat{\boldsymbol{\mu}}_j = A_j \mathbf{y} + \boldsymbol{b}_j$. For instance, the variance of $\hat{\boldsymbol{\mu}}_j$ is $\mathbb{E}\|\hat{\boldsymbol{\mu}}_j - \mathbb{E}\hat{\boldsymbol{\mu}}_j\|_2^2 = \sigma^2 \operatorname{Tr}(A_j^T A_j)$. If $\hat{\boldsymbol{\mu}}_j$ is a least squares estimator, $A_j$ is an orthoprojector, and the variance becomes $\sigma^2 \operatorname{Tr} A_j$. Thus, the statistical complexity of the estimator $\hat{\boldsymbol{\mu}}_j$ clearly appears in the right-hand side of the oracle inequality (1.3) proved in [10]. Thus, one may think that the price to pay for aggregating affine estimators, that is, the function $\mathrm{PRICE}_M(\delta)$ in (1.2), depends on the statistical complexity of the estimators to aggregate.

The bound (1.3) may lead to the conclusion that the price to pay for aggregation of affine estimators can be substantially larger than $\sigma^2 \log(M/\delta)$ which is the price for aggregating deterministic vectors. Indeed, the extra term $4\sigma^2 \operatorname{Tr}(A_j)$ may be large in common situation where the trace of some matrices $A_j$ is large. For example, if one aggregates the estimators $\hat{\boldsymbol{\mu}}_1 = \lambda_1 \mathbf{y}, \ldots, \hat{\boldsymbol{\mu}}_M = \lambda_M \mathbf{y}$, for some positive real numbers $\lambda_1, \ldots, \lambda_M$, then the term $4\sigma^2 \operatorname{Tr}(A_j)$ in the above oracle inequality is of order $\sigma^2 n \lambda_j$ for each $j = 1, \ldots, M$, which can be greater than the optimal rate $\sigma^2 \log M$. This term $4\sigma^2 \operatorname{Tr}(A_j)$ makes the oracle inequality (1.3) suitable only for scenarios where the matrices $A_j$ have small trace. But more importantly, the term $\sigma^2 \operatorname{Tr} A_j$ suggests that the price to pay for aggregating affine estimators increases with the statistical complexities of the estimators to aggregate.

The results discussed above rely on specific assumptions on the matrices $A_1, \ldots, A_M$ [10, 12, 27]. This raises a third question, although not as important as the two questions above:

3. Does the nature of the matrices $A_1, \ldots, A_M$ have an impact on the price to pay to aggregate these affine estimators? Is the price in (1.2) substantially smaller if the matrices are orthoprojectors, semipositive definite or symmetric?

The main contribution of the present paper is to answer the three questions raised above:

1. It is proved in Theorem 2.1 that a penalized procedure over the simplex satisfies the sharp oracle inequality (1.2) with $\operatorname{PRICE}_M(\delta) = c\sigma^2 \log(M/\delta)$ for some absolute constant $c > 0$. This price is of the same order as for the problem of aggregation of deterministic vectors. Thus, the dependence between the estimators and the data used to aggregate them induces no extra cost.

2. The form of the affine estimators to aggregate has no impact on the price to pay for aggregation. In particular, the sharp oracle inequalities of the present paper do not involve quantities dependent on $A_j$ such as $\sigma^2 \operatorname{Tr} A_j$.

3. The only assumption made on the matrices $A_1, \ldots, A_M$ is that $\max_{j \neq k} \|\|A_j - A_k\|\|_2$ is bounded from above, where $\|\| \cdot \|\|_2$ is the operator norm. All other assumptions on the matrices $A_1, \ldots, A_M$ can be dropped, in particular the matrices can be nonsymmetric and have negative eigenvalues.

The paper is organized as follows. In Section 1.1, we define the notation used throughout the paper. Section 2 defines a penalized procedure over the simplex and shows that it achieves sharp oracle inequalities in deviation for aggregation of affine estimators. The role of the penalty is studied in Section 3 and Section 4. Prior weights are considered in Section 5. Section 6 shows that the estimator is robust to variance misspecification and to non-Gaussianity of the noise. Some examples are given in Section 7. Section 8 is devoted to the proofs.

1.1. *Notation.* Let $\mathbf{f} = (f_1, \ldots, f_n)^T \in \mathbf{R}^n$ be an unknown regression vector. We observe $n$ random variables (1.1) where $\xi_1, \ldots, \xi_n$ are sub-Gaussian random variables, with $\mathbb{E}[\xi_i] = 0$ and $\mathbb{E}[\xi_i^2] = \sigma^2$. It can be rewritten in the vector form $\mathbf{y} = \mathbf{f} + \boldsymbol{\xi}$ where $\mathbf{y} = (Y_1, \ldots, Y_n)^T$, $\mathbf{f} = (f_1, \ldots, f_n)^T$ and $\boldsymbol{\xi} = (\xi_1, \ldots, \xi_n)^T$.

For any estimator $\hat{\boldsymbol{\mu}}$ of $\mathbf{f}$, we measure the quality of estimation of $\mathbf{f}$ with the loss $\|\hat{\boldsymbol{\mu}} - \mathbf{f}\|_2^2$, where $\|\cdot\|_2$ is the Euclidean norm in $\mathbf{R}^n$. Let $M \geq 2$. We consider $M$ affine estimators of the form

$$\hat{\boldsymbol{\mu}}_j = A_j \mathbf{y} + \boldsymbol{b}_j, \qquad j = 1, \ldots, M.$$

The matrices $A_1, \ldots, A_M$ and the vectors $\boldsymbol{b}_1, \ldots, \boldsymbol{b}_M \in \mathbf{R}^n$ are deterministic. Define the simplex in $\mathbf{R}^M$:

$$\Lambda^M = \left\{ \boldsymbol{\theta} \in \mathbf{R}^M, \sum_{j=1}^M \theta_j = 1, \forall j = 1 \cdots M, \theta_j \geq 0 \right\}.$$

For any $\boldsymbol{\theta} \in \Lambda^M$, let

$$A_{\boldsymbol{\theta}} = \sum_{j=1}^M \theta_j A_j, \qquad \boldsymbol{b}_{\boldsymbol{\theta}} = \sum_{j=1}^M \theta_j \boldsymbol{b}_j, \qquad \hat{\boldsymbol{\mu}}_{\boldsymbol{\theta}} = A_{\boldsymbol{\theta}} \mathbf{y} + \boldsymbol{b}_{\boldsymbol{\theta}}.$$

Let $\boldsymbol{e}_1, \ldots, \boldsymbol{e}_M$ be the vectors of the canonical basis in $\mathbf{R}^M$. Then $\hat{\boldsymbol{\mu}}_j = \hat{\boldsymbol{\mu}}_{\boldsymbol{e}_j}$ for all $j = 1, \ldots, M$.

An orthoprojector is an $n \times n$ matrix $P$ such that

(1.4) $$P = P^T = P^2.$$

Denote by $I_{n \times n}$ the $n \times n$-identity matrix. For any $n \times n$ real matrix $A = (a_{i,j})_{i,j=1,\ldots,n}$, define the operator norm of $A$, the Frobenius (or Hilbert–Schmidt) norm of $A$ and the nuclear norm of $A$ respectively by

$$|\!|\!| A |\!|\!|_2 = \sup_{x \neq 0} \frac{\|Ax\|_2}{\|x\|_2}, \qquad \|A\|_F = \sqrt{\sum_{i,j=1,\ldots,n} a_{i,j}^2}.$$

The following inequalities hold for any two squared matrices $M, M'$:

(1.5) $$|\!|\!| MM' |\!|\!|_2 \leq |\!|\!| M |\!|\!|_2 |\!|\!| M' |\!|\!|_2, \qquad \|MM'\|_F \leq |\!|\!| M |\!|\!|_2 \|M'\|_F.$$

Finally, denote by log the natural logarithm with $\log(e) = 1$.

**2. A penalized procedure on the simplex.** For any $\boldsymbol{\theta} \in \Lambda^M$, define

(2.1) $$C_p(\boldsymbol{\theta}) := \|\hat{\boldsymbol{\mu}}_{\boldsymbol{\theta}}\|_2^2 - 2\mathbf{y}^T \hat{\boldsymbol{\mu}}_{\boldsymbol{\theta}} + 2\sigma^2 \operatorname{Tr}(A_{\boldsymbol{\theta}}),$$

which is Mallows [28] $C_p$-criterion. Next, define

(2.2) $$H_{\mathrm{pen}}(\boldsymbol{\theta}) = C_p(\boldsymbol{\theta}) + \frac{1}{2} \operatorname{pen}(\boldsymbol{\theta}),$$

where

(2.3) $$\operatorname{pen}(\boldsymbol{\theta}) = \sum_{j=1}^M \theta_j \|\hat{\boldsymbol{\mu}}_{\boldsymbol{\theta}} - \hat{\boldsymbol{\mu}}_j\|_2^2.$$

We consider the estimator $\hat{\boldsymbol{\mu}}_{\hat{\boldsymbol{\theta}}_{\mathrm{pen}}}$ where

(2.4) $$\hat{\boldsymbol{\theta}}_{\mathrm{pen}} \in \operatorname*{argmin}_{\boldsymbol{\theta} \in \Lambda^M} H_{\mathrm{pen}}(\boldsymbol{\theta}).$$

The function $H_{\text{pen}}$ is quadratic and convex (cf. Lemma 8.2). Minimizing $H_{\text{pen}}$ over the simplex is a convex quadratic program for which efficient algorithms are available. The convexity of $H_{\text{pen}}$ also proves that $\hat{\boldsymbol{\theta}}_{\text{pen}}$ is well defined, although it may not be unique (e.g., if all $\hat{\boldsymbol{\mu}}_j$ are the same then $H_{\text{pen}}$ is constant on the simplex).

We now explain the meaning of the terms that appear in (2.2). If $\boldsymbol{\theta}$ is fixed, $C_p(\boldsymbol{\theta})$ is an unbiased estimate of the quantity

$$(2.5) \qquad R(\boldsymbol{\theta}) := \|\hat{\boldsymbol{\mu}}_{\boldsymbol{\theta}}\|_2^2 - 2\mathbf{f}^T \hat{\boldsymbol{\mu}}_{\boldsymbol{\theta}} = \|\hat{\boldsymbol{\mu}}_{\boldsymbol{\theta}} - \mathbf{f}\|_2^2 - \|\mathbf{f}\|_2^2,$$

which is the quantity of interest $\|\hat{\boldsymbol{\mu}}_{\boldsymbol{\theta}} - \mathbf{f}\|_2^2$ up to the additive constant $\|\mathbf{f}\|_2^2$.

The penalty (2.3) is borrowed from the $Q$-aggregation procedure, which is a powerful tool to derive sharp oracle inequalities in deviation when the loss is strongly convex [4, 11, 26, 31]. Since the estimators $\hat{\boldsymbol{\mu}}_1, \ldots, \hat{\boldsymbol{\mu}}_M$ depend on the data, the penalty (2.3) is data-driven, which is not the case if $\hat{\boldsymbol{\mu}}_1, \ldots, \hat{\boldsymbol{\mu}}_M$ are deterministic vectors as in [11]. In order to give some geometric insights on the penalty (2.3), let $c \in \mathbf{R}^n$ be a solution of $M$ linear equations $2c^T \hat{\boldsymbol{\mu}}_j = \|\hat{\boldsymbol{\mu}}_j\|_2^2$, $j = 1, \ldots, M$, and assume only in the rest of this paragraph that such a solution exists, even though this assumption cannot be fulfilled for $M > n$. Then

$$(2.6) \quad \text{pen}(\boldsymbol{\theta}) = \sum_{j=1}^M \theta_j \|\hat{\boldsymbol{\mu}}_j\|_2^2 - \|\hat{\boldsymbol{\mu}}_{\boldsymbol{\theta}}\|_2^2 = 2c^T \hat{\boldsymbol{\mu}}_{\boldsymbol{\theta}} - \|\hat{\boldsymbol{\mu}}_{\boldsymbol{\theta}}\|_2^2 = \|c\|_2^2 - \|\hat{\boldsymbol{\mu}}_{\boldsymbol{\theta}} - c\|_2^2.$$

We can write $\text{pen}(\boldsymbol{\theta}) = g(\hat{\boldsymbol{\mu}}_{\boldsymbol{\theta}})$ for some function $g$ defined on the convex hull of $\{\hat{\boldsymbol{\mu}}_1, \ldots, \hat{\boldsymbol{\mu}}_M\}$. Equation (2.6) shows that the level sets of the function $g$ are Euclidean balls centered at $c$. The function $g$ is nonnegative. It is minimal at the extreme points $\hat{\boldsymbol{\mu}}_1, \ldots, \hat{\boldsymbol{\mu}}_M$ since $g(\hat{\boldsymbol{\mu}}_j) = 0$ for all $j = 1, \ldots, M$ and $g$ is maximal at the projection of $c$ on the convex hull of $\{\hat{\boldsymbol{\mu}}_1, \ldots, \hat{\boldsymbol{\mu}}_M\}$. Intuitively, the penalty (2.3) pushes $\boldsymbol{\theta}$ away from the center of the simplex towards the vertices. Thus, the level sets of the function $\boldsymbol{\theta} \to \text{pen}(\boldsymbol{\theta})$ in $\mathbf{R}^M$ are ellipsoids centered at $\boldsymbol{\theta}_c$, where $\boldsymbol{\theta}_c$ is the unique point in $\mathbf{R}^M$ such that $\hat{\boldsymbol{\mu}}_{\boldsymbol{\theta}_c} = c$. If $M > n$ or if the vector $c$ is not well defined, the level sets of $\text{pen}(\cdot)$ are more intricate and cannot be described in such a simple way.

THEOREM 2.1 (Main result). *Let $M \geq 2$. For $j = 1, \ldots, M$, consider the affine estimators $\hat{\boldsymbol{\mu}}_j = A_j \mathbf{y} + \boldsymbol{b}_j$ and let*

$$(2.7) \qquad \phi := \max\left(1, \max_{j,k=1,\ldots,M: j \neq k} \frac{1}{2} \|\!|A_j - A_k|\!\|_2\right).$$

*Assume that the noise random variables $\xi_1, \ldots \xi_n$ are i.i.d. $\mathcal{N}(0, \sigma^2)$. Let $\hat{\boldsymbol{\theta}}_{\text{pen}}$ be the estimator defined in (2.4). Then for all $x > 0$ the estimator $\hat{\boldsymbol{\mu}}_{\hat{\boldsymbol{\theta}}_{\text{pen}}}$ satisfies with probability greater than $1 - \exp(-x)$,*

$$(2.8) \qquad \|\hat{\boldsymbol{\mu}}_{\hat{\boldsymbol{\theta}}_{\text{pen}}} - \mathbf{f}\|_2^2 \leq \min_{j=1,\ldots,M} \|\hat{\boldsymbol{\mu}}_j - \mathbf{f}\|_2^2 + 30\phi^2\sigma^2(x + 2\log M).$$

*Furthermore*,

$$(2.9) \quad \mathbb{E}[\|\hat{\boldsymbol{\mu}}_{\hat{\boldsymbol{\theta}}_{\text{pen}}} - \mathbf{f}\|_2^2] \le \mathbb{E}\Big[\min_{j=1,\ldots,M} \|\hat{\boldsymbol{\mu}}_j - \mathbf{f}\|_2^2\Big] + 60\phi^2 \sigma^2 \log(M).$$

The sharp oracle inequality in deviation given in [10] presents an additive term proportional to $\sigma^2 \text{Tr}(A_j)$, as in (1.3). An improvement of the present paper is the absence of this additive term which can be large for matrices $A_j$ with large trace. Our analysis shows that the quantities $\sigma^2 \text{Tr}(A_j)$ are not meaningful for the problem of aggregation of affine estimators, and Theorem 2.1 improves upon the earlier result of [10].

The quantity $\phi$ defined in (2.7) appears in the right-hand side of the oracle inequalities of Theorem 2.1. Cohen [9] established that estimators of the form $\hat{\boldsymbol{\mu}}_j = A_j \mathbf{y}$ with $\|\|A_j\|\|_2 > 1$ are inadmissible. Thus, if $\hat{\boldsymbol{\mu}}_1, \ldots, \hat{\boldsymbol{\mu}}_M$ are admissible affine estimators, then $\|\|A_j\|\|_2 \le 1$ and the quantity (2.7) is equal to 1.

We relax all assumptions on the matrices $A_1, \ldots, A_M$, for instance, they may be nonsymmetric and have negative eigenvalues. The above result shows that the restrictions on the matrices $A_1, \ldots, A_M$ introduced in [10, 12, 27] are not intrinsic to the problem of aggregation of affine estimators.

The next proposition shows that the bounds of Theorem 2.1 are optimal in a minimax sense. For any $\mathbf{f} \in \mathbf{R}^n$, we denote by $\mathbb{P}_{\mathbf{f}}$ the probability measure of the random variable $\mathbf{y} = \mathbf{f} + \boldsymbol{\xi}$. A lower bound for aggregation of deterministic vectors was proved in [32], Theorem 5.4 with $S = 1$. This lower bound implies the following result.

PROPOSITION 2.1. *There exist absolute constants $c^*, C^*, p^* > 0$ such that the following holds. For all $M, n \ge C^*$, there exist $\boldsymbol{b}_1, \ldots, \boldsymbol{b}_M \in \mathbf{R}^n$ and orthoprojectors $A_1, \ldots, A_M$ of rank one such that*

$$(2.10) \quad \inf_{\hat{\boldsymbol{\mu}}} \sup_{\mathbf{f} \in \mathbf{R}^n} \mathbb{P}_{\mathbf{f}}\Big(\|\hat{\boldsymbol{\mu}} - \mathbf{f}\|_2^2 - \min_{k=1,\ldots,M} \|\boldsymbol{b}_k - \mathbf{f}\|_2^2 \ge c^* \sigma^2 \log(M)\Big) \ge p^*,$$

$$(2.11) \quad \inf_{\hat{\boldsymbol{\mu}}} \sup_{\mathbf{f} \in \mathbf{R}^n} \mathbb{P}_{\mathbf{f}}\Big(\|\hat{\boldsymbol{\mu}} - \mathbf{f}\|_2^2 - \min_{k=1,\ldots,M} \|A_k \mathbf{y} - \mathbf{f}\|_2^2 \ge c^* \sigma^2 \log(M)\Big) \ge p^*,$$

*where the infima are taken over all estimators $\hat{\boldsymbol{\mu}}$.*

This implies that the bounds of Theorem 2.1 are rate minimax in terms of the aggregation price. The lower bound can be constructed either with a dictionary of deterministic vectors [cf. (2.10)], or with a dictionary of orthoprojectors of rank one [cf. (2.11)].

**3. The penalty (2.3) improves upon model selection based on $C_p$.** In order to explain the role of the penalty (2.3) for the problem of aggregation of affine estimators, consider first the standard empirical risk minimization scheme based on the $C_p$ criterion. Define $\hat{J}$ as

$$\hat{J} \in \operatorname*{argmin}_{j=1,\ldots,M} C_p(\boldsymbol{e}_j), \tag{3.1}$$

where $C_p(\cdot)$ is defined in (2.1). Using that $C_p(\boldsymbol{e}_{\hat{j}}) \leq C_p(\boldsymbol{e}_k)$ for all $k = 1, \ldots, M$ together with the definition of $C_p(\cdot)$ and $R(\cdot)$ given in (2.1) and (2.5), the following holds almost surely:

$$\|\hat{\boldsymbol{\mu}}_{\hat{j}} - \mathbf{f}\|_2^2 \leq \min_{k=1,\ldots,M} \|\hat{\boldsymbol{\mu}}_k - \mathbf{f}\|_2^2 + \max_{j,k=1,\ldots,M} \Delta_{jk}, \tag{3.2}$$

where $\Delta_{jk} := C_p(\boldsymbol{e}_k) - C_p(\boldsymbol{e}_j) - (R(\boldsymbol{e}_k) - R(\boldsymbol{e}_j))$. Thus, it is possible to prove an oracle inequality for the estimator $\hat{\boldsymbol{\mu}}_{\hat{j}}$ if we can control the quantities $\Delta_{jk}$ uniformly over all pairs $j, k = 1, \ldots, M$. These quantities can be rewritten as

$$\begin{aligned}\Delta_{jk} &= 2\boldsymbol{\xi}^T((A_j - A_k)\mathbf{f} + \boldsymbol{b}_j - \boldsymbol{b}_k) \\ &\quad + 2\big(\boldsymbol{\xi}^T(A_j - A_k)\boldsymbol{\xi} - \sigma^2 \operatorname{Tr}(A_j - A_k)\big).\end{aligned} \tag{3.3}$$

Two stochastic terms appear in $\Delta_{jk}$. The first is a centered Gaussian random variable with variance $4\sigma^2 \|(A_j - A_k)\mathbf{f} + \boldsymbol{b}_j - \boldsymbol{b}_k\|_2^2$. The second is a centered quadratic form in $\boldsymbol{\xi}$, and it can be shown that its variance is of order $\sigma^4 \|A_j - A_k\|_F^2$. This quadratic term is sometimes called a Gaussian chaos of order 2. The deviations of these two terms are governed by the following concentration inequalities. For any vector $\boldsymbol{v} \in \mathbf{R}^n$, a standard Gaussian tail bound gives

$$\mathbb{P}\big(\boldsymbol{v}^T \boldsymbol{\xi} > \sigma \|\boldsymbol{v}\|_2 \sqrt{2x}\big) \leq \exp(-x) \qquad \forall x > 0. \tag{3.4}$$

For the Gaussian chaos of order 2, the following is proved in [7], Example 2.12.

LEMMA 3.1. *Assume that $\boldsymbol{\xi} \sim \mathcal{N}(0, \sigma^2 I_{n \times n})$. For any squared matrix $B$ of size $n$,*

$$\mathbb{P}\big(\boldsymbol{\xi}^T B \boldsymbol{\xi} - \sigma^2 \operatorname{Tr} B > 2\sigma^2 \|B\|_F \sqrt{x} + 2\sigma^2 \|\!|\!| B |\!|\!|_2 x\big) \leq \exp(-x), \tag{3.5}$$

*where $\sigma^2 \operatorname{Tr} B = \mathbb{E}[\boldsymbol{\xi}^T B \boldsymbol{\xi}]$.*

We set $\boldsymbol{v} = 2((A_j - A_k)\mathbf{f} + \boldsymbol{b}_j - \boldsymbol{b}_k)$ and $B = 2(A_k - A_j)$ to study the deviations of the random variable $\Delta_{jk}$. If $\|\!|\!| A_j - A_k |\!|\!|_2$ is small, (3.4) and (3.5) yield that the deviations of $\Delta_{jk}$ are of order of the two quantities

$$\sigma \|(A_j - A_k)\mathbf{f} + \boldsymbol{b}_j - \boldsymbol{b}_k\|_2, \qquad \sigma^2 \|A_j - A_k\|_F, \tag{3.6}$$

that is, the standard deviations of the two terms in $\Delta_{jk}$. The concentration inequalities (3.4) and (3.5) are known to be tight [23], thus there is little hope to bound the deviations of $\Delta_{jk}$ independently of $\mathbf{f}$, $A_j$ and $A_k$ in order to prove a sharp oracle inequality. It is possible to refine the above analysis and to prove the following oracle inequality, though with a leading constant greater than 1.

PROPOSITION 3.1. *There exist absolute constants $c, C > 0$ such that the following holds. Let $0 < \varepsilon < c$ and let $\hat{J}$ be the estimator defined in (3.1). For all $x > 0$, the estimator $\hat{\boldsymbol{\mu}}_{\hat{j}}$ satisfies with probability greater than $1 - 2\exp(-x)$:*

$$\text{(3.7)} \quad \|\hat{\boldsymbol{\mu}}_{\hat{j}} - \mathbf{f}\|_2^2 \le (1+\varepsilon) \min_{k=1,\ldots,M} \|\hat{\boldsymbol{\mu}}_k - \mathbf{f}\|_2^2 + C\phi^2 \sigma^2 (x + 2\log M)/\varepsilon,$$

*where $\phi$ is defined in (2.7).*

The proof of Proposition 3.1 is given in Section 8.6. The estimator $\hat{\boldsymbol{\mu}}_{\hat{j}}$ fails to achieve an oracle inequality with leading constant 1 [the leading constant in (3.7) is $1 + \varepsilon$] and with an error term of order $\sigma^2 \log M$. This drawback cannot be repaired for all procedures of the form $\hat{\boldsymbol{\mu}}_{\hat{K}}$ where $\hat{K}$ is an estimator valued in $\{1, \ldots, M\}$. Indeed, it is proved in [16], Section 6.4.2 and Proposition 6.1, that there exist $\mathbf{f}_1, \mathbf{f}_2 \in \mathbf{R}^n$ and orthoprojectors $A_1, A_2$ such that for any estimator $\hat{K}$ valued in $\{1, 2\}$

$$\sup_{\mathbf{f} \in \{\mathbf{f}_1, \mathbf{f}_2\}} \left( \mathbb{E} \|A_{\hat{K}} \mathbf{y} - \mathbf{f}\|_2^2 - \min_{j=1,2} \mathbb{E}\|A_j \mathbf{y} - \mathbf{f}\|_2^2 \right) \ge \sigma^2 \sqrt{n}/4,$$

provided that $n$ is larger than some absolute constant. Inspection of the proof of this result reveals that

$$\sigma \|(A_2 - A_1)\mathbf{f} + \boldsymbol{b}_2 - \boldsymbol{b}_1\|_2 \ge \sigma^2 \sqrt{n} \qquad \forall \mathbf{f} \in \{\mathbf{f}_1, \mathbf{f}_2\},$$

where we set $\boldsymbol{b}_1 = \boldsymbol{b}_2 = 0$. Thus, this lower bound of order $\sqrt{n}$ is related to the Gaussian component of the random variable $\Delta_{12}$, that is, to the term $\boldsymbol{\xi}^T((A_1 - A_2)\mathbf{f} + \boldsymbol{b}_1 - \boldsymbol{b}_2)$; cf. (3.3).

The procedure $\hat{\boldsymbol{\mu}}_{\hat{j}}$ fails to achieve a sharp oracle inequality because the variances of the two components of $\Delta_{jk}$ may be large and cannot be controlled. The role of the penalty (2.3) is exactly to control the deviations of $\Delta_{jk}$ by controlling the terms (3.6). The following proposition makes this precise.

PROPOSITION 3.2. *Let $\hat{\boldsymbol{\theta}}_{\text{pen}}$ be the estimator (2.4). Then almost surely,*

$$\text{(3.8)} \quad \|\hat{\boldsymbol{\mu}}_{\hat{\boldsymbol{\theta}}_{\text{pen}}} - \mathbf{f}\|_2^2 \le \min_{q=1,\ldots,M} (\|\hat{\boldsymbol{\mu}}_q - \mathbf{f}\|_2^2) + \max_{j,k=1,\ldots,M} \left( \Delta_{jk} - \frac{1}{2} \|\hat{\boldsymbol{\mu}}_j - \hat{\boldsymbol{\mu}}_k\|_2^2 \right),$$

*where $\Delta_{jk}$ is the quantity (3.3). Furthermore, for all $j, k = 1, \ldots, M$,*

$$\text{(3.9)} \quad \mathbb{E}\left[\frac{1}{2}\|\hat{\boldsymbol{\mu}}_j - \hat{\boldsymbol{\mu}}_k\|_2^2\right] = \frac{1}{2}\|(A_j - A_k)\mathbf{f} + \boldsymbol{b}_j - \boldsymbol{b}_k\|_2^2 + \frac{\sigma^2}{2}\|A_j - A_k\|_F^2.$$

The proof of (3.8) is given in Section 4 below. A bias-variance decomposition directly yields (3.9), since $\mathbb{E}[\hat{\boldsymbol{\mu}}_j - \hat{\boldsymbol{\mu}}_k] = (A_j - A_k)\mathbf{f} + \boldsymbol{b}_j - \boldsymbol{b}_k$ and $\mathbb{E}\|\hat{\boldsymbol{\mu}}_j - \hat{\boldsymbol{\mu}}_k - \mathbb{E}[\hat{\boldsymbol{\mu}}_j - \hat{\boldsymbol{\mu}}_k]\|_2^2 = \mathbb{E}\|(A_j - A_k)\boldsymbol{\xi}\|_2^2 = \sigma^2 \|A_j - A_k\|_F^2$.

Compared with (3.2), the right-hand side of (3.8) presents the quantities $-\frac{1}{2}\|\hat{\boldsymbol{\mu}}_j - \hat{\boldsymbol{\mu}}_k\|_2^2$. We will explain below that these quantities appear because of the interplay between the penalty (2.3) and the strong convexity of $H_{\text{pen}}$.

From (3.8), an outline of the proof of Theorem 2.1 is as follows. By combining the simple inequality (8.5) and Proposition 8.1 below, we will prove that for any pair $(j, k)$ we have

$$\mathbb{E} \exp\left(\lambda_0 \left(\Delta_{jk} - \frac{1}{2}\|\hat{\boldsymbol{\mu}}_j - \hat{\boldsymbol{\mu}}_k\|_2^2\right)\right) \leq 1$$

for $\lambda_0 = (30\phi^2\sigma^2)^{-1}$ if the noise $\boldsymbol{\xi}$ has distribution $\mathcal{N}(0, \sigma^2 I_{n\times n})$. Thus, one has

$$\mathbb{E} \exp\left(\lambda_0 \max_{j,k=1,\ldots,M}\left(\Delta_{jk} - \frac{1}{2}\|\hat{\boldsymbol{\mu}}_j - \hat{\boldsymbol{\mu}}_k\|_2^2\right)\right) \leq M^2.$$

Then Jensen's inequality yields (2.9) while a Chernoff bound yields (2.8). This explains the success of the penalty (2.3) for the problem of model selection type aggregation: the penalty and the strong convexity of $H_{\text{pen}}$ provide the quantity $-\frac{1}{2}\|\hat{\boldsymbol{\mu}}_j - \hat{\boldsymbol{\mu}}_k\|_2^2$, and this quantity is exactly what is needed to control the deviations of the random variable $\Delta_{jk}$.

**4. Strong convexity and the penalty** (2.3). To further understand the interplay between the penalty (2.3) and the strong convexity of $H_{\text{pen}}$, we now give the proof of (3.8).

PROOF OF (3.8). Let $k = 1, \ldots, M$ be fixed. The simplex $\Lambda^M$ is a convex set and the function $H_{\text{pen}}$ is convex, hence we have

$$\nabla H_{\text{pen}}(\hat{\boldsymbol{\theta}}_{\text{pen}})^T(\boldsymbol{e}_k - \hat{\boldsymbol{\theta}}_{\text{pen}}) \geq 0;$$

cf. [8], Section 4.2.3, equation (4.21). Inequality (3.8) follows from

$$\|\hat{\boldsymbol{\mu}}_{\hat{\boldsymbol{\theta}}_{\text{pen}}} - \mathbf{f}\|_2^2 - \|\hat{\boldsymbol{\mu}}_k - \mathbf{f}\|_2^2$$

(4.1) $$\leq \|\hat{\boldsymbol{\mu}}_{\hat{\boldsymbol{\theta}}_{\text{pen}}} - \mathbf{f}\|_2^2 - \|\hat{\boldsymbol{\mu}}_k - \mathbf{f}\|_2^2 + \nabla H_{\text{pen}}(\hat{\boldsymbol{\theta}}_{\text{pen}})^T(\boldsymbol{e}_k - \hat{\boldsymbol{\theta}}_{\text{pen}}),$$

(4.2) $$= \sum_{j=1}^M \hat{\theta}_{\text{pen},j}\left(\Delta_{jk} - \frac{1}{2}\|\hat{\boldsymbol{\mu}}_j - \hat{\boldsymbol{\mu}}_k\|_2^2\right),$$

(4.3) $$\leq \max_{j=1,\ldots,M}\left(\Delta_{jk} - \frac{1}{2}\|\hat{\boldsymbol{\mu}}_j - \hat{\boldsymbol{\mu}}_k\|_2^2\right).$$

Equality (4.2) is obtained by simple algebra while (4.3) is a consequence of $\sum_{j=1}^M \hat{\theta}_{\text{pen},j} = 1$ and $\hat{\theta}_{\text{pen},j} \geq 0$ for all $j = 1, \ldots, M$. □

It is possible to interpret this argument in light of the interplay between strong convexity and the penalty (2.3). The right-hand side of (4.1) satisfies

$$\|\hat{\boldsymbol{\mu}}_{\hat{\boldsymbol{\theta}}_{\text{pen}}} - \mathbf{f}\|_2^2 - \|\hat{\boldsymbol{\mu}}_k - \mathbf{f}\|_2^2 + \nabla H_{\text{pen}}(\hat{\boldsymbol{\theta}}_{\text{pen}})^T (\boldsymbol{e}_k - \hat{\boldsymbol{\theta}}_{\text{pen}})$$

$$= \sum_{j=1}^{M} \hat{\theta}_{\text{pen},j} \Delta_{jk} - \frac{1}{2}[\text{pen}(\hat{\boldsymbol{\theta}}_{\text{pen}}) + \|\hat{\boldsymbol{\mu}}_{\hat{\boldsymbol{\theta}}_{\text{pen}}} - \hat{\boldsymbol{\mu}}_k\|_2^2].$$

The term $\|\hat{\boldsymbol{\mu}}_{\hat{\boldsymbol{\theta}}_{\text{pen}}} - \hat{\boldsymbol{\mu}}_k\|_2^2$ comes from the strong convexity of the function $H_{\text{pen}}$. By simple algebra or using (8.10) with $\boldsymbol{g} = \hat{\boldsymbol{\mu}}_k$, we have

(4.4) $\quad \text{pen}(\hat{\boldsymbol{\theta}}_{\text{pen}}) + \underbrace{\|\hat{\boldsymbol{\mu}}_{\hat{\boldsymbol{\theta}}_{\text{pen}}} - \hat{\boldsymbol{\mu}}_k\|_2^2}_{\text{Term given by the strong convexity of } H_{\text{pen}}} = \sum_{j=1}^{M} \hat{\theta}_{\text{pen},j} \underbrace{\|\hat{\boldsymbol{\mu}}_j - \hat{\boldsymbol{\mu}}_k\|_2^2}_{\text{Term that controls the deviations of } \Delta_{jk}}.$

Formula (4.4) highlights a feature of the penalty (2.3): the penalty transforms the quadratic term given by strong convexity into the linear term given by the right-hand side of (4.4).

The strong convexity of $C_p(\cdot)$ and $H_{\text{pen}}(\cdot)$ is understood with respect to the pseudometric

$$\|\hat{\boldsymbol{\mu}}_{\boldsymbol{\theta}} - \hat{\boldsymbol{\mu}}_{\boldsymbol{\theta}'}\|_2, \qquad \boldsymbol{\theta}, \boldsymbol{\theta}' \in \mathbf{R}^M,$$

so it is not the strong convexity in the Euclidean norm. We say that a function $V(\cdot)$ is strongly convex with coefficient $\gamma > 0$ over the simplex if for all $\boldsymbol{\theta}, \boldsymbol{\theta}' \in \Lambda^M$,

$$V(\boldsymbol{\theta}) \geq V(\boldsymbol{\theta}') + \nabla V(\boldsymbol{\theta}')^T (\boldsymbol{\theta} - \boldsymbol{\theta}') + \gamma \|\hat{\boldsymbol{\mu}}_{\boldsymbol{\theta}} - \hat{\boldsymbol{\mu}}_{\boldsymbol{\theta}'}\|_2^2.$$

The strong convexity of $H_{\text{pen}}$ could be used because $H_{\text{pen}}$ is minimized over the simplex and not just over the vertices. Indeed, minimizing a strongly convex function over a discrete set, as in the definition of $\hat{J}$, only grants the inequalities

$$C_p(\boldsymbol{e}_{\hat{j}}) \leq C_p(\boldsymbol{e}_k) \qquad \text{for all } k = 1, \ldots, M.$$

Because the simplex is a convex set, minimizing the strongly convex function $H_{\text{pen}}$ over the simplex grants the inequalities

$$H_{\text{pen}}(\hat{\boldsymbol{\theta}}_{\text{pen}}) \leq H_{\text{pen}}(\boldsymbol{\theta}) - \frac{1}{2}\|\hat{\boldsymbol{\mu}}_{\boldsymbol{\theta}} - \hat{\boldsymbol{\mu}}_{\hat{\boldsymbol{\theta}}_{\text{pen}}}\|_2^2 \qquad \text{for all } \boldsymbol{\theta} \in \Lambda^M.$$

One could also consider the estimator $\hat{\boldsymbol{\theta}}_C \in \operatorname{argmin}_{\boldsymbol{\theta} \in \Lambda^M} C_p(\boldsymbol{\theta})$. Because of the strong convexity of $C_p(\cdot)$, this estimator enjoys the inequalities:

$$C_p(\hat{\boldsymbol{\theta}}_C) \leq C_p(\boldsymbol{\theta}) - \|\hat{\boldsymbol{\mu}}_{\boldsymbol{\theta}} - \hat{\boldsymbol{\mu}}_{\hat{\boldsymbol{\theta}}_{\text{pen}}}\|_2^2 \qquad \text{for all } \boldsymbol{\theta} \in \Lambda^M.$$

The above displays highlight the fact that $C_p(\cdot)$ and $H_{\text{pen}}(\cdot)$ have different strong convexity coefficients. This is because $H_{\text{pen}}(\cdot) = C_p(\cdot) + (1/2)\text{pen}(\cdot)$ and $(1/2)\text{pen}(\cdot)$ is strongly concave with coefficient $1/2$, thus the strong convexity coefficient of $H_{\text{pen}}(\cdot)$ is less than that of $C_p(\cdot)$. We refer to Lemma 8.2 for a rigorous proof of the strong convexity of $H_{\text{pen}}$ and $C_p$.

The estimator $\hat{\boldsymbol{\mu}}_{\hat{\boldsymbol{\theta}}_C}$ is another candidate for the problem of aggregation of affine estimators. It is close to the estimator $\hat{\boldsymbol{\mu}}_{\hat{\boldsymbol{\theta}}_{\text{pen}}}$, except that the penalty (2.3) has been removed from the function to minimize. It was proved in [11], Section 2.2, that the estimator $\hat{\boldsymbol{\mu}}_{\hat{\boldsymbol{\theta}}_C}$ performs poorly: for large enough $M$ and $n$, there exist $\mathbf{f}$ and $\boldsymbol{b}_1, \ldots, \boldsymbol{b}_M \in \mathbf{R}^n$ such that with probability greater than $1/4$,

$$\|\hat{\boldsymbol{\mu}}_{\hat{\boldsymbol{\theta}}_C} - \mathbf{f}\|_2^2 \geq \min_{j=1,\ldots,M} \|\hat{\boldsymbol{\mu}}_j - \mathbf{f}\|_2^2 + \frac{\sigma^2 \sqrt{n}}{48},$$

where $\hat{\boldsymbol{\mu}}_j = \boldsymbol{b}_j$ for all $j = 1, \ldots, M$.

**5. Prior weights.** We consider now the problem of aggregation of $M$ affine estimators with a prior probability distribution $\boldsymbol{\pi} = (\pi_1, \ldots, \pi_M)^T$ on the finite set of indices $\{1, \ldots, M\}$.

THEOREM 5.1. *Let $M \geq 2$. For $j = 1, \ldots, M$, consider the estimator $\hat{\boldsymbol{\mu}}_j = A_j \mathbf{y} + \boldsymbol{b}_j$ and let $\phi$ be defined in (2.7). Let $\boldsymbol{\pi} = (\pi_1, \ldots, \pi_M)^T \in \Lambda^M$. Assume that the noise $\boldsymbol{\xi}$ has distribution $\mathcal{N}(0, \sigma^2 I_{n \times n})$. Let $\hat{\boldsymbol{\theta}}_{\boldsymbol{\pi}} \in \arg\min_{\boldsymbol{\theta} \in \Lambda^M} V_{\text{pen}}(\boldsymbol{\theta})$ where*

$$(5.1) \qquad V_{\text{pen}}(\boldsymbol{\theta}) := H_{\text{pen}}(\boldsymbol{\theta}) + 30\phi^2 \sigma^2 \sum_{j=1}^M \theta_j \log \frac{1}{\pi_j}.$$

*Then for all $x > 0$, with probability greater than $1 - \exp(-x)$,*

$$(5.2) \qquad \|\hat{\boldsymbol{\mu}}_{\hat{\boldsymbol{\theta}}_{\boldsymbol{\pi}}} - \mathbf{f}\|_2^2 \leq \min_{j=1,\ldots,M} \left(\|\hat{\boldsymbol{\mu}}_j - \mathbf{f}\|_2^2 + 60\phi^2 \sigma^2 \log \frac{1}{\pi_j}\right) + 30\phi^2 \sigma^2 x.$$

*Furthermore,*

$$(5.3) \qquad \mathbb{E}\|\hat{\boldsymbol{\mu}}_{\hat{\boldsymbol{\theta}}_{\boldsymbol{\pi}}} - \mathbf{f}\|_2^2 \leq \mathbb{E} \min_{j=1,\ldots,M} \left(\|\hat{\boldsymbol{\mu}}_j - \mathbf{f}\|_2^2 + 60\phi^2 \sigma^2 \log \frac{1}{\pi_j}\right).$$

The prior probability distribution $\boldsymbol{\pi} = (\pi_j)_{j=1,\ldots,M}$ is deterministic and does not depend on the data $\mathbf{y} = (Y_1, \ldots, Y_n)^T$. The only difference between the function (2.2) and the function minimized in (5.1) is the term proportional to

$$(5.4) \qquad \phi^2 \sigma^2 \sum_{j=1}^M \theta_j \log \frac{1}{\pi_j}.$$

This term allows us to weight the candidates $\hat{\boldsymbol{\mu}}_1, \ldots, \hat{\boldsymbol{\mu}}_M$ with the prior probability distribution $(\pi_j)_{j=1,\ldots,M}$ based on some prior knowledge about the estimators $\hat{\boldsymbol{\mu}}_1, \ldots, \hat{\boldsymbol{\mu}}_M$. For example, if the estimators are orthoprojectors, one can set prior weights that decrease with the rank the orthoprojectors [32, 34]. The same term is used in [26] whereas [10] uses the Kullback–Leibler divergence of $\boldsymbol{\theta}$ from $\boldsymbol{\pi}$. It is shown in [11] that for aggregation of deterministic vectors, one may use a quantity of the form $\sum_{j=1}^{M} \theta_j \log(\rho(\theta_j)/\pi_j)$ where $\rho(\cdot)$ satisfies $\rho(t) \geq t$ and $t \to t \log(\rho(t))$ is convex. This suggests that we could use the Kullback–Leibler divergence of $\boldsymbol{\theta}$ from $\boldsymbol{\pi}$ instead of (5.4), but in their current form our proofs only hold with the "linear entropy" (5.4).

**6. Robustness of the estimator $\hat{\boldsymbol{\mu}}_{\hat{\boldsymbol{\theta}}_{\text{pen}}}$.** We prove in this section that the procedure (2.4) is robust to non-Gaussian noise distributions and to variance misspecification.

6.1. *Robustness to non-Gaussian noise.* The following result shows that the penalized procedure (2.4) is robust to non-Gaussian noise distributions.

THEOREM 6.1. *Let $M \geq 2$. Let $\bar{\sigma} > 0$. For $j = 1, \ldots, M$, consider the estimator $\hat{\boldsymbol{\mu}}_j = A_j \mathbf{y} + \boldsymbol{b}_j$ and assume that $\|\|A_j\|\|_2 \leq 1$ for all $j = 1, \ldots, M$. Assume that the noise components $\xi_1, \ldots, \xi_n$ are i.i.d., centered with variance $\sigma^2$ and satisfy for all $\boldsymbol{b} \in \mathbf{R}^n$, all matrices $B$ and all $x > 0$:*

$$\mathbb{P}(\boldsymbol{\xi}^T \boldsymbol{b} > \bar{\sigma}\sqrt{2x}) \leq \exp(-x), \tag{6.1}$$

$$\mathbb{P}(\boldsymbol{\xi}^T B \boldsymbol{\xi} - \sigma^2 \operatorname{Tr} B > 2\sigma\bar{\sigma}\|B\|_F\sqrt{x} + 2\bar{\sigma}^2\|\|B\|\|_2 x) \leq \exp(-x). \tag{6.2}$$

*Let $\hat{\boldsymbol{\theta}}_{\text{pen}}$ be the estimator defined in (2.4). Then for all $x > 0$, the estimator $\hat{\boldsymbol{\mu}}_{\hat{\boldsymbol{\theta}}_{\text{pen}}}$ satisfies with probability greater than $1 - 2\exp(-x)$,*

$$\|\hat{\boldsymbol{\mu}}_{\hat{\boldsymbol{\theta}}_{\text{pen}}} - \mathbf{f}\|_2^2 \leq \min_{j=1,\ldots,M} \|\hat{\boldsymbol{\mu}}_j - \mathbf{f}\|_2^2 + 46\bar{\sigma}^2(2\log M + x). \tag{6.3}$$

Let $K > 0$. If the random variables $\xi_1, \ldots, \xi_n$ are i.i.d., centered with variance $\sigma^2$ and $K$-sub-Gaussian in the sense that $\log \mathbb{E}[e^{t\xi_i}] \leq K^2 t^2/2$ for all $t \in \mathbf{R}$ and all $i = 1, \ldots, n$, then (6.1) is satisfied with $\bar{\sigma} = cK$ for some absolute constant $c > 0$ [40], Section 5.2.3. As $\sigma \leq K$, (6.1) is also satisfied with $\bar{\sigma} = cK^2/\sigma$. By the Hanson–Wright inequality [20, 35, 41], (6.2) also holds with $\bar{\sigma} = cK^2/\sigma$ for another absolute constant $c > 0$. Thus, for i.i.d. $K$-sub-Gaussian random variables with variance $\sigma^2$, (6.3) yields

$$\|\hat{\boldsymbol{\mu}}_{\hat{\boldsymbol{\theta}}_{\text{pen}}} - \mathbf{f}\|_2^2 \leq \min_{j=1,\ldots,M} \|\hat{\boldsymbol{\mu}}_j - \mathbf{f}\|_2^2 + C(K^4/\sigma^2)(2\log M + x), \tag{6.4}$$

for some absolute constant $C > 0$. For most common examples of sub-Gaussian random variables, the standard deviation $\sigma$ is of the same order as the sub-Gaussian norm $K$, so the bound (6.4) is satisfying. This bound may not be tight if the standard deviation is pathologically small compared to the sub-Gaussian norm.

6.2. *Robustness to variance misspecification.* In order to construct the estimator (2.4) by minimizing (2.2), the knowledge of the variance of the noise is needed. However, the following proposition shows that the procedure (2.4) is robust to variance misspecification, that is, the result still holds if the variance is replaced by an estimator $\hat\sigma^2$ as soon as $\hat\sigma^2$ is consistent in a weak sense defined below.

THEOREM 6.2 (Aggregation under variance misspecification). *Let $M \geq 2$. For $j = 1, \ldots, M$, consider the estimator $\hat{\boldsymbol{\mu}}_j = A_j \mathbf{y} + \boldsymbol{b}_j$. Assume that the noise random variables $\xi_1, \ldots \xi_n$ are i.i.d. $\mathcal{N}(0, \sigma^2)$. Let $\hat\sigma^2$ be an estimator and assume that*

(6.5) $\quad \forall j = 1, \ldots, M, \qquad A_j = A_j^T = A_j^2, \qquad \delta := \mathbb{P}(|\sigma^2 - \hat\sigma^2| > \sigma^2/8) < 1.$

*Let $\hat{\boldsymbol{\theta}}_{\hat\sigma} = \mathrm{argmin}_{\boldsymbol{\theta} \in \Lambda^M} W_{\mathrm{pen}}(\boldsymbol{\theta})$ where*

(6.6) $\qquad W_{\mathrm{pen}}(\boldsymbol{\theta}) := \|\hat{\boldsymbol{\mu}}_{\boldsymbol{\theta}}\|_2^2 - 2\mathbf{y}^T \hat{\boldsymbol{\mu}}_{\boldsymbol{\theta}} + 2\hat\sigma^2 \mathrm{Tr}(A_{\boldsymbol{\theta}}) + \frac{1}{2}\mathrm{pen}(\boldsymbol{\theta}).$

*Then for all $x > 0$, with probability greater than $1 - \delta - \exp(-x)$,*

$$\|\hat{\boldsymbol{\mu}}_{\hat{\boldsymbol{\theta}}_{\hat\sigma}} - \mathbf{f}\|_2^2 \leq \min_{j=1,\ldots,M} \|\hat{\boldsymbol{\mu}}_j - \mathbf{f}\|_2^2 + 48\sigma^2(x + 2\log M).$$

The proof of Theorem 6.2 is given in Section 8.3. In (6.5), the matrices $A_1, \ldots, A_M$ are assumed to be orthoprojectors, so Theorem 6.2 is a result for aggregation of least squares estimators. As soon as an estimator $\hat\sigma^2$ satisfies with high probability $|\hat\sigma^2 - \sigma^2| \leq \sigma^2/8$, optimal aggregation of least squares estimators is possible. This condition is weaker than consistency, as any estimator $\hat\sigma^2$ that converges to $\sigma^2$ in probability satisfies this condition for $n$ large enough.

The proof of Theorem 6.2 exploits the form of the penalty (2.3) and the strong convexity of the function (6.6). Similar to Proposition 3.2, we will prove that almost surely

(6.7)
$$\|\hat{\boldsymbol{\mu}}_{\hat{\boldsymbol{\theta}}_{\hat\sigma}} - \mathbf{f}\|_2^2 \leq \min_{q=1,\ldots,M} \|\hat{\boldsymbol{\mu}}_q - \mathbf{f}\|_2^2$$
$$+ \max_{j,k=1,\ldots,M} \Big(\Delta_{jk} - \frac{1}{2}\|\hat{\boldsymbol{\mu}}_j - \hat{\boldsymbol{\mu}}_k\|_2^2$$
$$+ 2(\sigma^2 - \hat\sigma^2)\mathrm{Tr}(A_j - A_k)\Big),$$

where $\Delta_{jk}$ is the quantity (3.3). The only difference from (3.8) is in the extra term $2(\sigma^2 - \hat\sigma^2)\mathrm{Tr}(A_j - A_k)$ that appears because we used $\hat\sigma^2$ instead of $\sigma^2$ in the definition of $W_{\mathrm{pen}}(\cdot)$. On the event $|\hat\sigma^2 - \sigma^2| \leq \sigma^2/8$, it is easy to check that (cf. Lemma 8.1)

$$2(\sigma^2 - \hat\sigma^2)\mathrm{Tr}(A_j - A_k) \leq \frac{\sigma^2}{4}\|A_j - A_k\|_F^2.$$

As explained in the discussion that follows Proposition 3.2, the quantity $\frac{1}{2}\|\hat{\boldsymbol{\mu}}_j - \hat{\boldsymbol{\mu}}_k\|_2^2$ is given by the interplay between the penalty (2.3) and the strong convexity of the function that is minimized. By (3.9), the expectation of this quantity is greater than $(\sigma^2/2)\|A_j - A_k\|_F^2$. Thus, the penalty (2.3) and the strong convexity of $W_{\text{pen}}$ provide exactly what is needed to compensate the difference between $\hat{\sigma}^2$ and $\sigma^2$. Hence, the proof of Theorem 6.2 reveals that the robustness to variance misspecification is in fact due to the interplay between the penalty (2.3) and the strong convexity of $W_{\text{pen}}$.

The papers [3, 17, 18] aim at performing aggregation of least squares estimators when $\sigma^2$ is unknown, but unlike Theorem 6.2 the oracle inequalities that they established have a leading constant greater than 1. To our knowledge, Theorem 6.2 is the first aggregation result, with leading constant 1, that is robust to variance misspecification.

In the following, we describe several situations where the suitable estimator $\hat{\sigma}^2$ is available.

EXAMPLE 6.1 (An estimator $\hat{\sigma}^2$ that does not depend on **y**). In [12], Section 3.1, two contexts are given where an unbiased estimator of the covariance matrix, independent from **y**, is available. For example, the noise level can be estimated independently if the signal is captured multiple times by a single device, or if several identical devices capture the same signal.

EXAMPLE 6.2 (Difference based estimators). In nonparametric regression where the nonrandom design points are equi-spaced in [0, 1], a well-known estimator of the noise level is the difference based estimator $1/(2n-2)\sum_{i=1}^{n-1}(y_{i+1}-y_i)^2$. This technique can be refined with more complex difference sequences [13, 19], and extends to design points in a multidimensional space [29]. For images, where the underlying space is 2-dimensional, there exist efficient methods which require no multiplication [21].

EXAMPLE 6.3 (Consistent estimation of $\sigma^2$ in high-dimensional linear regression). In a high-dimensional setting, it is possible to estimate $\sigma^2$ under classical assumptions in high-dimensional regression. First, the scaled Lasso [36] allows a joint estimation of the regression coefficients and of the noise level $\sigma^2$. The estimator $\hat{\sigma}^2$ of the scaled Lasso converges in probability to the true noise level $\sigma^2$ [36], Theorem 1, and $\hat{\sigma}^2/\sigma^2$ is asymptotically normal [36], equation (19). Second, [5] proposes to estimate $\sigma^2$ with with a recursive procedure that uses Lasso residuals, and nonasymptotic guarantees are proved in the supplementary material [5]. Third, [6] provides nonasymptotic bounds on the estimation of $\sigma^2$ by the residuals of the square-root Lasso ([6], Theorem 2), and these bounds imply consistency. In Theorem 6.2, we require that $|\hat{\sigma}^2/\sigma^2 - 1| \leq 1/8$ with high-probability and this requirement is far weaker than the guarantees obtained in [5, 36].

## 7. Examples.

7.1. *Adaptation to the smoothness.* For all $n \geq 1$, given continuous parameters $\beta \geq 1$ and $L > 0$, we consider subsets $\Theta(\beta, L) \subset \mathbf{R}^n$. We assume that for each $\beta \geq 1$, there exists a squared matrix $A_\beta$ of size $n$ with $|||A_\beta|||_2 \leq 1$ such that for all $L > 0$, as $n \to +\infty$,

$$(7.1) \quad \inf_{\hat{\mathbf{f}}} \sup_{\mathbf{f} \in \Theta(\beta, L)} \frac{1}{n} \mathbb{E}\|\mathbf{f} - \hat{\mathbf{f}}\|_2^2 \sim \sup_{\mathbf{f} \in \Theta(\beta, L)} \frac{1}{n} \mathbb{E}\|\mathbf{f} - A_\beta \mathbf{y}\|_2^2 \sim C^* n^{\frac{-2\beta}{2\beta+1}},$$

where $a_n \sim b_n$ if and only if $a_n/b_n \to 1$ as $n \to +\infty$, the infimum is taken over all estimators and the constant $C^* > 0$ may depend on $\beta, L$ and $\sigma$. The above assumption holds for Sobolev ellipsoids in nonparametric regression, and in this case one can choose the Pinsker filters for the matrices $A_\beta$ (cf. [38], Theorem 3.2). For Sobolev ellipsoids, there exist different estimators that adapt to the unknown smoothness [12, 15, 38].

Consider the following aggregation procedure. Assume that $n \geq 3$ and let $M = \lceil 120 \log(n)(\log \log n)^2 \rceil$. For all $j = 1, \ldots, M$, let

$$\beta_j = (1 + 1/(\log(n) \log \log n))^{j-1}.$$

We aggregate the linear estimators $(\hat{\boldsymbol{\mu}}_j = A_{\beta_j} \mathbf{y})_{j=1,\ldots,M}$ using the procedure (2.4) of Theorem 2.1, and denote by $\tilde{\boldsymbol{\mu}}$ the resulting estimator. The following adaptation result is a direct consequence of Theorem 2.1.

PROPOSITION 7.1. *For all $n \geq 3$, $\beta \geq 1$ and $L > 0$, let $\Theta(\beta, L) \subset \mathbf{R}^n$ such that as $n \to +\infty$, (7.1) is satisfied for some matrices $A_\beta$ with $|||A_\beta|||_2 \leq 1$. Assume that the sets $\Theta(\beta, L)$ are ordered, that is, $\Theta(\beta, L) \subset \Theta(\beta', L)$ for any $\beta > \beta'$ and any $L > 0$. For all $\beta \geq 1$ and $L > 0$, the estimator $\tilde{\boldsymbol{\mu}}$ defined above satisfies as $n \to +\infty$*

$$\lim_{n \to +\infty} \sup_{\mathbf{f} \in \Theta(\beta, L)} \frac{1}{n} \mathbb{E}\|\mathbf{f} - \tilde{\boldsymbol{\mu}}\|_2^2 n^{\frac{2\beta}{2\beta+1}} = C^*.$$

Proposition 7.1 is proved in Section 8.7. The above procedure adapts to the unknown smoothness in exact asymptotic sense by aggregating only

$$\log(n)(\log \log n)^2$$

estimators so its computational complexity is small. Another feature is that the minimax rate and the minimax constant $C^*$ are not altered by the aggregation step.

7.2. *The best convex combination as a benchmark.* We consider convex combinations of the estimators $\hat{\boldsymbol{\mu}}_1, \ldots, \hat{\boldsymbol{\mu}}_M$ to construct the estimator (2.4). The goal of this section is to study the performance of the estimator (2.4) if the benchmark is $\min_{\boldsymbol{\theta} \in \Lambda^M} \|\hat{\boldsymbol{\mu}}_{\boldsymbol{\theta}} - \mathbf{f}\|_2^2$ instead of $\min_{k=1,\ldots,M} \|\hat{\boldsymbol{\mu}}_k - \mathbf{f}\|_2^2$.

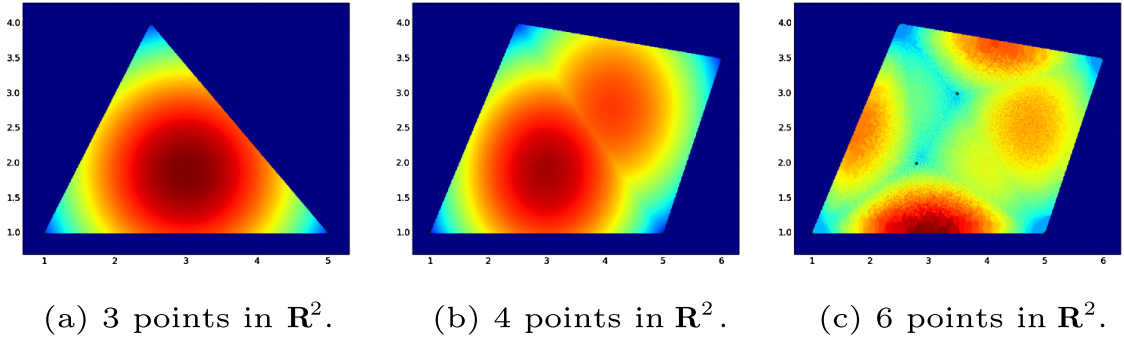

(a) 3 points in $\mathbf{R}^2$.  (b) 4 points in $\mathbf{R}^2$.  (c) 6 points in $\mathbf{R}^2$.

FIG. 1. *Penalty* (2.3) *heatmaps. Largest penalty in red, smallest in blue.*

The penalty (2.3) vanishes at the extreme points: $\text{pen}(e_j) = 0$ for all $j = 1, \ldots, M$, and it pushes $\hat{\boldsymbol{\mu}}_{\hat{\boldsymbol{\theta}}_{\text{pen}}}$ towards the points $\{\hat{\boldsymbol{\mu}}_1, \ldots, \hat{\boldsymbol{\mu}}_M\}$. This can be seen in Figure 1. Consider a noise-free problem where $\sigma = 0$. Let $\mathbf{f} \in \mathbf{R}^n$ and let $\boldsymbol{v} \in \mathbf{R}^n$ be such that $\boldsymbol{v}^T \mathbf{f} = 0$ and $\|\boldsymbol{v}\|_2^2 = 3\|\mathbf{f}\|_2^2 > 0$. Let also $\rho = 4\|\boldsymbol{v}\|_2^2$. Consider estimators $\hat{\boldsymbol{\mu}}_1 = 2\mathbf{f}$, $\hat{\boldsymbol{\mu}}_2 = \mathbf{f} + \boldsymbol{v}$ and $\hat{\boldsymbol{\mu}}_3 = \mathbf{f} - \boldsymbol{v}$. These estimators are such that $\|\hat{\boldsymbol{\mu}}_j\|_2^2 = \rho > 0$ for all $j = 1, \ldots, 3$ (here $M =$ and the estimators are deterministic because $\sigma = 0$). Then by simple algebra we have $\text{pen}(\boldsymbol{\theta}) = \rho - \|\hat{\boldsymbol{\mu}}_{\boldsymbol{\theta}}\|_2^2$ and $H_{\text{pen}}(\boldsymbol{\theta}) = (1/2)\|\hat{\boldsymbol{\mu}}_{\boldsymbol{\theta}} - 2\mathbf{f}\|_2^2 + c$ where $c$ is constant that depends on $\mathbf{f}$ but not on $\boldsymbol{\theta}$. For this example, although $\mathbf{f}$ lies in the convex hull of $\{\hat{\boldsymbol{\mu}}_1, \hat{\boldsymbol{\mu}}_2, \hat{\boldsymbol{\mu}}_3\}$, the estimator $\hat{\boldsymbol{\mu}}_{\hat{\boldsymbol{\theta}}}$ defined in (2.4) will be equal to $2\mathbf{f}$ instead of $\mathbf{f}$ and is likely to be a bad procedure with respect to the benchmark $\min_{\boldsymbol{\theta} \in \Lambda^M} \|\hat{\boldsymbol{\mu}}_{\boldsymbol{\theta}} - \mathbf{f}\|_2^2$. This fact is not surprising since the penalty penalizes heavily some regions of the convex hull of the estimators. Furthermore, this procedure is tailored for the benchmark $\min_{k=1,\ldots,M} \|\hat{\boldsymbol{\mu}}_k - \mathbf{f}\|_2^2$ and its goal is not to mimic the best convex combination of the estimators.

It is possible to modify the procedure (2.4) to construct an estimator that performs well with respect to the best convex combination of $M$ linear estimators. Let

$$(7.2) \qquad m := \left\lfloor \sqrt{\frac{n}{\log(1 + M/\sqrt{n})}} \right\rfloor.$$

If $m \geq 1$, define the set $\Lambda_m^M \subset \Lambda^M$ as

$$(7.3) \qquad \Lambda_m^M := \left\{ \frac{1}{m} \sum_{q=1}^m \boldsymbol{u}_q, \boldsymbol{u}_1, \ldots, \boldsymbol{u}_m \in \{\boldsymbol{e}_1, \ldots, \boldsymbol{e}_M\} \right\}.$$

Denote by $|\Lambda_m^M|$ the cardinality of $\Lambda_m^M$. We aggregate the affine estimators $(\hat{\boldsymbol{\mu}}_{\boldsymbol{u}})_{\boldsymbol{u} \in \Lambda_m^M}$ using the procedure (2.4) and denote by $\hat{\boldsymbol{\mu}}_{\Lambda_m^M}$ the resulting estimator.

PROPOSITION 7.2. *Let $M, n \geq 1$. For $j = 1, \ldots, M$, consider the estimator $\hat{\boldsymbol{\mu}}_j = A_j \mathbf{y} + \boldsymbol{b}_j$ for any $n \times n$ matrix $A_j$ and vector $\boldsymbol{b}_j \in \mathbf{R}^n$. Assume that $\boldsymbol{\xi} \sim$*

$\mathcal{N}(0, \sigma^2 I_{n \times n})$ and that for some constant $R > 0$,

$$\frac{1}{n}\|\mathbf{f}\|_2^2 \leq R^2, \qquad \frac{1}{n}\|\mathbf{b}_j\|_2^2 \leq R^2, \qquad \|\|A_j\|\|_2 \leq 1 \qquad \forall j = 1, \ldots, M.$$

For all $x > 0$, the estimator $\hat{\boldsymbol{\theta}}_C \in \operatorname{argmin}_{\boldsymbol{\theta} \in \Lambda^M} C_p(\boldsymbol{\theta})$ satisfies with probability greater than $1 - 2\exp(-x)$,

$$
(7.4) \quad \frac{1}{n}\|\hat{\boldsymbol{\mu}}_{\hat{\boldsymbol{\theta}}_C} - \mathbf{f}\|_2^2 \leq \min_{\boldsymbol{\theta} \in \Lambda^M} \frac{1}{n}\|\hat{\boldsymbol{\mu}}_{\boldsymbol{\theta}} - \mathbf{f}\|_2^2 + 8(\sigma^2 + \sigma R\sqrt{2})\sqrt{\frac{x + 2\log M}{n}} + \frac{8\sigma^2(x + 2\log M)}{n}.
$$

If $M \leq \sqrt{n}(\exp(n) - 1)$, then for all $x > 0$, the estimator $\hat{\boldsymbol{\mu}}_{\Lambda_m^M}$ defined above satisfies with probability greater than $1 - 3\exp(-x)$,

$$
(7.5) \quad \frac{1}{n}\|\hat{\boldsymbol{\mu}}_{\Lambda_m^M} - \mathbf{f}\|_2^2 \leq \min_{\boldsymbol{\theta} \in \Lambda^M} \frac{1}{n}\|\hat{\boldsymbol{\mu}}_{\boldsymbol{\theta}} - \mathbf{f}\|_2^2 + C\max(R^2, \sigma^2)\sqrt{\frac{\log(1 + M/\sqrt{n})}{n}} + \frac{C\sigma^2 x}{n}.
$$

Proposition 7.2 is proved in Section 8.8. To our knowledge, this is the first result that provides a sharp oracle inequality for the problem of aggregation of affine estimators with respect to the convex oracle. However, there is a large literature on convex aggregation when the estimators to aggregate are deterministic, which corresponds to the particular case $A_j = 0$ for all $j = 1, \ldots, M$. When the error is measured with the scaled squared norm $\frac{1}{n}\|\cdot\|_2^2$, the minimax rate of convex aggregation is known to be of order $M/n$ if $M \leq \sqrt{n}$ and $\sqrt{\log(1 + M/\sqrt{n})/n}$ if $M > \sqrt{n}$. For our setting, this is proved in [32]. This elbow effect was first established for regression with random design [37] and then extended to other settings in [31, 33]. All these results assume that the estimators to aggregate are deterministic or independent of the data used for aggregation. The lower bound [32], Theorem 5.3 with $S = M$, $\delta = \sigma$ and $R = \log(1 + eM)$, yields that there exist absolute constants $c, C > 0$ such that if $\log(1 + eM)^2 \leq Cn$, there exist deterministic vectors $\hat{\boldsymbol{\mu}}_1 = \mathbf{b}_1, \ldots, \hat{\boldsymbol{\mu}}_M = \mathbf{b}_M$ such that for all estimators $\hat{\boldsymbol{\mu}}$,

$$\sup_{\mathbf{f} \in \mathbf{R}^n} \mathbb{P}_{\mathbf{f}}\left(\frac{1}{n}\|\hat{\boldsymbol{\mu}} - \mathbf{f}\|_2^2 - \min_{\boldsymbol{\theta} \in \Lambda^M}\frac{1}{n}\|\hat{\boldsymbol{\mu}}_{\boldsymbol{\theta}} - \mathbf{f}\|_2^2 \geq c\sigma^2\left(\frac{M}{n} \wedge \sqrt{\frac{\log(1 + M/\sqrt{n})}{n}}\right)\right) \geq c.$$

Thus, if $M \geq \sqrt{n}$, (7.5) is optimal in a minimax sense up to absolute constants, and (7.4) is optimal up to logarithmic factors. However, we do not know whether the minimax rate is $M/n$ when $M < \sqrt{n}$, as in the case of aggregation of deterministic vectors.

The problem of linear aggregation of affine estimators remains open. It is only known that for linear aggregation of deterministic vectors, the least squares estimator on a linear space of dimension $M$ achieves the rate $\sigma^2 M/n$, which is optimal in a minimax sense [31–33, 37].

A bound similar to (7.4) can be obtained by the Lasso estimator for for aggregation of deterministic vectors, that is, in the case $A_j = 0$ and $\hat{\boldsymbol{\mu}}_j = \boldsymbol{b}_j$ for all $j = 1, \ldots, M$. Let $\boldsymbol{x}_j = (\sqrt{n}/\|\boldsymbol{b}_j\|_2)\boldsymbol{b}_j$ for all $j = 1, \ldots, M$ and let $\mathbb{X}$ be the $n \times M$ matrix with columns $\boldsymbol{x}_1, \ldots, \boldsymbol{x}_M$. The Lasso estimator is defined by

$$\hat{\boldsymbol{\beta}} \in \operatorname*{argmin}_{\boldsymbol{\beta} \in \mathbf{R}^M} \frac{1}{2n}\|\mathbf{y} - \mathbb{X}\boldsymbol{\beta}\|_2^2 + 4\sigma\sqrt{\log(M)/n}\|\boldsymbol{\beta}\|_1,$$

where $\|\boldsymbol{\beta}\|_1 = \sum_{j=1}^M |\beta_j|$. This estimator satisfies with probability at least $1 - 1/M$ the oracle inequality

$$\frac{1}{n}\|\mathbf{f} - \mathbb{X}\hat{\boldsymbol{\beta}}\|_2^2 \le \min_{\boldsymbol{\beta} \in \mathbf{R}^M} \left(\frac{1}{n}\|\mathbf{f} - \mathbb{X}\boldsymbol{\beta}\|_2^2 + c\sigma\|\boldsymbol{\beta}\|_1\sqrt{\log(M)/n}\right),$$

where $c > 0$ is a numerical constant. This oracle inequality is proved in Theorem 4 in Sun and Zhang [36]. If we assume $\frac{1}{n}\|\boldsymbol{b}_j\|^2 \le R^2$ for all $j = 1, \ldots, M$ as in Proposition 7.2 above, then the right-hand side of the previous display is bounded from above by

$$\min_{\boldsymbol{\theta} \in \Lambda^M} \left(\frac{1}{n}\|\mathbf{f} - \boldsymbol{b}_{\boldsymbol{\theta}}\|_2^2 + c\sigma R\sqrt{\log(M)/n}\right).$$

Hence, the Lasso estimator satisfies a bound similar to the oracle inequalities of Proposition 7.2 for aggregation of deterministic vectors, that is, if $A_j = 0$ and $\hat{\boldsymbol{\mu}}_j = \boldsymbol{b}_j$ for all $j = 1, \ldots, M$.

## 8. Proofs.

8.1. *Preliminaries.* The following notation will be useful. Define for all $j, k = 1, \ldots, M$:

$$(8.1) \quad Q_{j,k} := \left(-2I_{n \times n} - \frac{1}{2}(A_k - A_j)^T\right)(A_k - A_j),$$

$$(8.2) \quad \boldsymbol{v}_{j,k} := (-2I_{n \times n} - (A_k - A_j)^T)((A_k - A_j)\mathbf{f} + \boldsymbol{b}_k - \boldsymbol{b}_j).$$

Let $B_{jk} = A_k - A_j$, so that $\hat{\boldsymbol{\mu}}_k - \hat{\boldsymbol{\mu}}_j = B_{jk}\boldsymbol{\xi} + (B_{jk}\mathbf{f} + \boldsymbol{b}_k - \boldsymbol{b}_j)$. Then

$$\|\hat{\boldsymbol{\mu}}_k - \hat{\boldsymbol{\mu}}_j\|_2^2 = \|B_{jk}\boldsymbol{\xi}\|_2^2 + \|B_{jk}\mathbf{f} + \boldsymbol{b}_k - \boldsymbol{b}_j\|_2^2 + 2\boldsymbol{\xi}^T B_{jk}^T(B_{jk}\mathbf{f} + \boldsymbol{b}_k - \boldsymbol{b}_j).$$

Thus, simple algebra yields that the quantity $\Delta_{jk}$ defined in (3.3) satisfies

$$\begin{aligned}\Delta_{jk} - \frac{1}{2}\|\hat{\boldsymbol{\mu}}_k - \hat{\boldsymbol{\mu}}_j\|_2^2 &= \boldsymbol{\xi}^T Q_{j,k}\boldsymbol{\xi} - \mathbb{E}[\boldsymbol{\xi}^T Q_{j,k}\boldsymbol{\xi}] + \boldsymbol{\xi}^T \boldsymbol{v}_{j,k}\\ &\quad - \frac{\sigma^2}{2}\|A_j - A_k\|_F^2 - \frac{1}{2}\|(A_k - A_j)\mathbf{f} + \boldsymbol{b}_k - \boldsymbol{b}_j\|_2^2,\end{aligned}$$

(8.3)

where we used the equality $\sigma^2 \|A_j - A_k\|_F^2 = \mathbb{E}[\|(A_j - A_k)\boldsymbol{\xi}\|_2^2]$ and the above definitions of $Q_{j,k}$ and $\boldsymbol{v}_{j,k}$. Furthermore, using (1.5) and (2.7) we have

(8.4)
$$\|\|Q_{j,k}\|\|_2 \leq 6\phi^2, \qquad \|Q_{j,k}\|_F \leq 3\phi \|A_k - A_j\|_F,$$
$$\|\boldsymbol{v}_{j,k}\|_2 \leq 4\phi \|(A_k - A_j)\mathbf{f} + \boldsymbol{b}_k - \boldsymbol{b}_j\|_2$$

for all $j, k = 1, \ldots, M$. This yields that

(8.5)
$$\Delta_{jk} - \frac{1}{2}\|\hat{\boldsymbol{\mu}}_k - \hat{\boldsymbol{\mu}}_j\|_2^2 \leq \boldsymbol{\xi}^T Q_{j,k} \boldsymbol{\xi} - \mathbb{E}[\boldsymbol{\xi}^T Q_{j,k} \boldsymbol{\xi}] + \boldsymbol{\xi}^T \boldsymbol{v}_{j,k}$$
$$- \frac{\sigma^2}{18\phi^2}\|Q_{j,k}\|_F^2 - \frac{1}{32\phi^2}\|\boldsymbol{v}_{j,k}\|_2^2.$$

PROPOSITION 8.1. *Let $\boldsymbol{v} \in \mathbf{R}^n$ and let $Q$ be any squared matrix of size $n$. Assume that $\xi_1, \ldots, \xi_n$ are i.i.d. $\mathcal{N}(0, \sigma^2)$ random variables. Then for all $u > 0$ such that $2u\sigma^2 \|\|Q\|\|_2 < 1$ we have*

(8.6)
$$\mathbb{E}\left[e^{u(\boldsymbol{\xi}^T Q \boldsymbol{\xi} - \mathbb{E}[\boldsymbol{\xi}^T Q \boldsymbol{\xi}] + \boldsymbol{\xi}^T \boldsymbol{v})}\right] \leq \exp\left(u^2 \sigma^2 \left(\frac{\sigma^2 \|Q\|_F^2 + \frac{\|\boldsymbol{v}\|_2^2}{2}}{1 - 2\sigma^2 \|\|Q\|\|_2 u}\right)\right).$$

*Furthermore, for some $\phi \geq 1$, define*

$$Z_{Q,\boldsymbol{v}} := \boldsymbol{\xi}^T Q \boldsymbol{\xi} - \mathbb{E}[\boldsymbol{\xi}^T Q \boldsymbol{\xi}] + \boldsymbol{\xi}^T \boldsymbol{v} - \frac{\sigma^2}{18\phi^2}\|Q\|_F^2 - \frac{1}{32\phi^2}\|\boldsymbol{v}\|_2^2,$$

$$Y_{Q,\boldsymbol{v}} := \boldsymbol{\xi}^T Q \boldsymbol{\xi} - \mathbb{E}[\boldsymbol{\xi}^T Q \boldsymbol{\xi}] + \boldsymbol{\xi}^T \boldsymbol{v} - \frac{\sigma^2}{36\phi^2}\|Q\|_F^2 - \frac{1}{32\phi^2}\|\boldsymbol{v}\|_2^2.$$

*If $\|\|Q\|\|_2 \leq 6\phi^2$, then for $u = 1/(30\phi^2\sigma^2)$ and $u' = 1/(48\phi^2\sigma^2)$ we have*

$$\mathbb{E}[e^{uZ_{Q,\boldsymbol{v}}}] \leq 1, \qquad \mathbb{E}[e^{u'Y_{Q,\boldsymbol{v}}}] \leq 1.$$

The proof relies on an argument similar to that of [24], Lemma 1.

PROOF. If $Q$ is not symmetric, let $Q_s = (Q + Q^T)/2$. We have $\|Q_s\|_F \leq \|Q\|_F$, $\|\|Q_s\|\|_2 \leq \|\|Q\|\|_2$ and almost surely $\boldsymbol{\xi}^T Q \boldsymbol{\xi} = \boldsymbol{\xi}^T Q_s \boldsymbol{\xi}$ so that if (8.6) holds for $Q_s$ then

$$\mathbb{E}\left[e^{u(\boldsymbol{\xi}^T Q \boldsymbol{\xi} - \mathbb{E}[\boldsymbol{\xi}^T Q \boldsymbol{\xi}] + \boldsymbol{\xi}^T \boldsymbol{v})}\right] \leq e^{u^2\sigma^2(\frac{\sigma^2\|Q_s\|_F^2 + \frac{\|\boldsymbol{v}\|_2^2}{2}}{1-2\sigma^2\|\|Q_s\|\|_2 u})} \leq e^{u^2\sigma^2(\frac{\sigma^2\|Q\|_F^2 + \frac{\|\boldsymbol{v}\|_2^2}{2}}{1-2\sigma^2\|\|Q\|\|_2 u})}.$$

Thus, the result for the symmetric matrix $Q_s$ implies the result for $Q$.

We now assume that $Q$ is symmetric. There exists a matrix $P$ with $P^T P = PP^T = I_{n\times n}$ such that $Q = P^T \operatorname{diag}(\lambda_1, \ldots, \lambda_n)P$, where $\lambda_1, \ldots, \lambda_n$ are the eigenvalues of $Q$. Let $\boldsymbol{w} = (1/\sigma)P\boldsymbol{v}$ and define the random variables $g_1, \ldots, g_n$ by $(g_1, \ldots, g_n)^T = (1/\sigma) P \boldsymbol{\xi}$. By the rotational invariance of the Gaussian distribution, $g_1, \ldots, g_n$ are i.i.d. $\mathcal{N}(0, 1)$ random variables. Thus, the random variable $\boldsymbol{\xi}^T Q \boldsymbol{\xi} - \mathbb{E}[\boldsymbol{\xi}^T Q \boldsymbol{\xi}] + \boldsymbol{\xi}^T \boldsymbol{v}$ has the same distribution as

$$\sigma^2 \sum_{i=1}^n W_i \qquad \text{where } W_i := \lambda_i (g_i^2 - 1) + g_i w_i.$$

For all $i = 1, \ldots, n$ and for all $t > 0$ such that $\max_{i=1,\ldots,n} 2t|\lambda_i| < 1$, integration using the probability density function of $g_i$ yields

$$\mathbb{E}[e^{t W_i}] = \frac{1}{\sqrt{1 - 2\lambda_i t}} e^{\frac{t^2 w_i^2}{2(1 - 2\lambda_i t)} - t\lambda_i} \leq e^{\frac{\lambda_i^2 t^2}{1 - 2|\lambda_i| t} + \frac{t^2 w_i^2}{2(1 - 2\lambda_i t)}},$$

where we used the inequalities

$$\log\left(\frac{1}{\sqrt{1 - 2v}}\right) \leq v + \frac{v^2}{1 - 2v} = v + \frac{v^2}{1 - 2|v|} \qquad \text{for all } v \in [0, 1/2),$$

$$\log\left(\frac{1}{\sqrt{1 - 2v}}\right) \leq v + v^2 \leq v + \frac{v^2}{1 - 2|v|} \qquad \text{for all } v \in (-1/2, 0].$$

This can be shown by comparing the power series expansions. As $|\lambda_i| \leq \|\!|Q|\!\|_2$ for all $i = 1, \ldots, n$, by independence of $W_1, \ldots, W_n$ we obtain

$$\mathbb{E}[e^{t \sum_{i=1}^n W_i}] \leq \exp\left(t^2 \left(\frac{\|Q\|_F^2 + \frac{\|\boldsymbol{w}\|_2^2}{2}}{1 - 2\|\!|Q|\!\|_2 t}\right)\right).$$

By definition of $\boldsymbol{w}$, we have $\|\boldsymbol{v}\|_2 = \sigma \|\boldsymbol{w}\|_2$, so setting $t = u\sigma^2$ completes the proof of (8.6).

The claims about $Z_{Q,\boldsymbol{v}}$ and $Y_{Q,\boldsymbol{v}}$ are direct consequences of (8.6). $\square$

8.2. *Proof of the main results.*

PROOF OF THEOREM 2.1. By Proposition 3.2, it is enough to prove that

(8.7) $$D_1 := \mathbb{E}[e^{u \max_{j,k=1,\ldots,M} (\Delta_{jk} - \frac{1}{2} \|\hat{\boldsymbol{\mu}}_j - \hat{\boldsymbol{\mu}}_k\|_2^2)}] \leq M^2$$

for $u = 1/(30\phi^2 \sigma^2)$. Then, Jensen's inequality yields (2.9) and a Chernoff bound yields (2.8).

We now prove (8.7). By (8.4), for all $j, k = 1, \ldots, M$ we have $\|\!|Q_{j,k}|\!\|_2 \leq 6\phi^2$. Using (8.5) and Proposition 8.1, we have

$$D_1 \leq \sum_{j=1}^M \sum_{k=1}^M \mathbb{E}[e^{u(\Delta_{jk} - \frac{1}{2} \|\hat{\boldsymbol{\mu}}_j - \hat{\boldsymbol{\mu}}_k\|_2^2)}] \leq \sum_{j=1}^M \sum_{k=1}^M \mathbb{E}[e^{u Z_{Q_{j,k}, \boldsymbol{v}_{j,k}}}] \leq M^2,$$

where for any matrix $Q$ and any $\boldsymbol{v} \in \mathbf{R}^n$, the random variable $Z_{Q,\boldsymbol{v}}$ is defined in Proposition 8.1. $\square$

PROOF OF THEOREM 5.1. Let $\beta = 30\phi^2\sigma^2$. Let $\hat{\boldsymbol{\theta}} = \hat{\boldsymbol{\theta}}_\pi$ for notational simplicity. The only difference between $H_{\text{pen}}$ and $V_{\text{pen}}$ is the linear term (5.4). As in the proof of (3.8) in Section 4, by convexity of $V_{\text{pen}}$ we have that for all $k = 1, \ldots, M$,

$$\|\hat{\boldsymbol{\mu}}_{\hat{\boldsymbol{\theta}}} - \mathbf{f}\|_2^2 - \|\hat{\boldsymbol{\mu}}_k - \mathbf{f}\|_2^2$$
$$\leq \|\hat{\boldsymbol{\mu}}_{\hat{\boldsymbol{\theta}}} - \mathbf{f}\|_2^2 - \|\hat{\boldsymbol{\mu}}_k - \mathbf{f}\|_2^2 + \nabla V_{\text{pen}}(\hat{\boldsymbol{\theta}})^T (\boldsymbol{e}_k - \hat{\boldsymbol{\theta}})$$
$$= 2\beta \log \frac{1}{\pi_k} + \sum_{j=1}^{M} \hat{\theta}_j \left( \Delta_{jk} - \frac{1}{2}\|\hat{\boldsymbol{\mu}}_j - \hat{\boldsymbol{\mu}}_k\|_2^2 - \beta \log \frac{1}{\pi_j \pi_k} \right)$$
$$\leq 2\beta \log \frac{1}{\pi_k} + \max_{j=1,\ldots,M} \left( \Delta_{jk} - \frac{1}{2}\|\hat{\boldsymbol{\mu}}_j - \hat{\boldsymbol{\mu}}_k\|_2^2 - \beta \log \frac{1}{\pi_j \pi_k} \right),$$

where $\Delta_{jk}$ is defined in (3.3). For all $u > 0$, let

$$D_2 := \mathbb{E}\left[ \exp\left( u \max_{j,k=1,\ldots,M} \left( \Delta_{jk} - \frac{1}{2}\|\hat{\boldsymbol{\mu}}_j - \hat{\boldsymbol{\mu}}_k\|_2^2 - \beta \log \frac{1}{\pi_j \pi_k} \right) \right) \right].$$

We now bound from above this moment generating function using (8.5) and Proposition 8.1. If $u = 1/\beta = 1/(30\phi^2\sigma^2)$, then

(8.8)
$$D_2 \leq \sum_{j=1}^{M} \sum_{k=1}^{M} \pi_j \pi_k \mathbb{E}\left[ e^{u(\Delta_{jk} - \frac{1}{2}\|\hat{\boldsymbol{\mu}}_j - \hat{\boldsymbol{\mu}}_k\|_2^2)} \right]$$
$$\leq \sum_{j=1}^{M} \sum_{k=1}^{M} \pi_j \pi_k \mathbb{E}\left[ e^{u Z_{Q_{j,k}, \boldsymbol{v}_{j,k}}} \right] \leq \sum_{j=1}^{M} \sum_{k=1}^{M} \pi_j \pi_k = 1.$$

As in the proof of Theorem 2.1, Jensen's inequality yields (5.2) while a Chernoff bound completes the proof of (5.2). □

PROOF OF THEOREM 6.1. For a fixed pair $(j,k)$, we apply (6.1) to the vector $\boldsymbol{v}_{j,k}$ and (6.2) to the matrix $Q_{j,k}$. Using (8.4),

$$\boldsymbol{\xi}^T Q_{j,k} \boldsymbol{\xi} - \mathbb{E}[\boldsymbol{\xi}^T Q_{j,k} \boldsymbol{\xi}] \leq \bar{\sigma}^2 12 x + 6\sigma\bar{\sigma}\|A_k - A_j\|_F \sqrt{x}$$
$$\leq 30\bar{\sigma}^2 x + \frac{\sigma^2}{2}\|A_k - A_j\|_F^2,$$
$$\boldsymbol{\xi}^T \boldsymbol{v}_{j,k} \leq \bar{\sigma} 4 \|(A_k - A_j)\mathbf{f} + \boldsymbol{b}_k - \boldsymbol{b}_j\|_2 \sqrt{2x}$$
$$\leq 16\bar{\sigma}^2 x + \frac{1}{2}\|(A_k - A_j)\mathbf{f} + \boldsymbol{b}_k - \boldsymbol{b}_j\|_2^2.$$

Combining this bound with (6.7), (8.3) and the union bound completes the proof. □

8.3. *Proof of Theorem* 6.2. The following inequality will be useful.

LEMMA 8.1 (Projection matrices). *Let $A$, $B$ be two squared matrices of size $n$ with $A^T = A = A^2$ and $B^T = B = B^2$. Then*

(8.9) $$|\operatorname{Tr}(A - B)| \le \|A - B\|_F^2.$$

PROOF. Without loss of generality, assume that $\operatorname{Tr} A \ge \operatorname{Tr} B$. As $\|A - B\|_F^2 = \|A\|_F^2 + \|B\|_F^2 - 2\operatorname{Tr}(AB)$ and $\|A\|_F^2 = \operatorname{Tr} A$, (8.9) is equivalent to $2\operatorname{Tr}(AB) \le 2\operatorname{Tr}(B)$. Notice that for projection matrices, $\operatorname{Tr}(AB) = \|AB\|_F^2 \le \|A\|_2^2 \|B\|_F^2 \le \|B\|_F^2 = \operatorname{Tr}(B)$ and the proof is complete. □

PROOF OF THEOREM 6.2. Let $\hat{\boldsymbol{\theta}} = \hat{\boldsymbol{\theta}}_{\hat{\sigma}}$ for notational simplicity. Since the matrices $A_j$ are orthoprojectors, the quantity $\phi$ defined in (2.7) is equal to 1. As in the proof of (3.8) in Section 4, by convexity of $W_{\text{pen}}$ we have that for all $k = 1, \ldots, M$,

$$\|\hat{\boldsymbol{\mu}}_{\hat{\boldsymbol{\theta}}} - \mathbf{f}\|_2^2 - \|\hat{\boldsymbol{\mu}}_k - \mathbf{f}\|_2^2$$
$$\le \|\hat{\boldsymbol{\mu}}_{\hat{\boldsymbol{\theta}}} - \mathbf{f}\|_2^2 - \|\hat{\boldsymbol{\mu}}_k - \mathbf{f}\|_2^2 + \nabla W_{\text{pen}}(\hat{\boldsymbol{\theta}})^T (\mathbf{e}_k - \hat{\boldsymbol{\theta}})$$
$$= \sum_{j=1}^M \hat{\theta}_j \left( \Delta_{jk} + 2(\hat{\sigma}^2 - \sigma^2)\operatorname{Tr}(A_j - A_k) - \frac{1}{2}\|\hat{\boldsymbol{\mu}}_j - \hat{\boldsymbol{\mu}}_k\|_2^2 \right)$$
$$\le \max_{j=1,\ldots,M} \left( \Delta_{jk} + 2(\hat{\sigma}^2 - \sigma^2)\operatorname{Tr}(A_j - A_k) - \frac{1}{2}\|\hat{\boldsymbol{\mu}}_j - \hat{\boldsymbol{\mu}}_k\|_2^2 \right) =: D_3,$$

where $\Delta_{jk}$ is defined in (3.3). The assumption on $\hat{\sigma}^2$ and (8.9) yield that on an event $\Omega_0$ of probability greater than $1 - \delta$,

$$2|(\hat{\sigma}^2 - \sigma^2)\operatorname{Tr}(A_j - A_k)| \le \frac{\sigma^2}{4}\|A_j - A_k\|_F^2 \qquad \text{for all } j, k = 1, \ldots, M.$$

Using (8.3) and (8.4), we obtain that on the event $\Omega_0$,

$$D_3 \le \max_{j,k=1,\ldots,M} \left( \boldsymbol{\xi}^T Q_{j,k} \boldsymbol{\xi} - \mathbb{E}[\boldsymbol{\xi}^T Q_{j,k} \boldsymbol{\xi}] + \boldsymbol{\xi}^T \boldsymbol{v}_{j,k} - \frac{\sigma^2}{36}\|Q_{j,k}\|_F^2 - \frac{1}{32}\|\boldsymbol{v}_{j,k}\|_2^2 \right)$$
$$= \max_{j,k=1,\ldots,M} Y_{Q_{j,k}, \boldsymbol{v}_{j,k}},$$

where $Q_{j,k}$ and $\boldsymbol{v}_{j,k}$ are defined in (8.1) and (8.2) while $Y_{Q,\boldsymbol{v}}$ is defined in Proposition 8.1 for any matrix $Q$ and any $\boldsymbol{v} \in \mathbf{R}^n$. Using Proposition 8.1, for $u = 1/(48\sigma^2)$ we have

$$\mathbb{E}\left[ \exp\left( u \max_{j,k=1,\ldots,M} Y_{Q_{j,k}, \boldsymbol{v}_{j,k}} \right) \right] \le \sum_{j=1}^M \sum_{k=1}^M \mathbb{E}[\exp(u Y_{Q_{j,k}, \boldsymbol{v}_{j,k}})] \le M^2.$$

By a Chernoff bound, this proves that on an event $\Omega_1$ of probability greater than $1 - e^{-x}$, we have $\max_{j,k=1,\ldots,M} Y_{Q_{j,k}, v_{j,k}},  \le 48\sigma^2(x + 2\log M)$. On the event $\Omega_0 \cap \Omega_1$, we have $D_3 \le 48\sigma^2(x + 2\log M)$ and the union bound yields that $\mathbb{P}(\Omega_0 \cap \Omega_1) \ge 1 - e^{-x} - \delta$. □

8.4. *Strong convexity.* The penalty (2.3) satisfies for any $g \in \mathbf{R}^n$ and any $\boldsymbol{\theta} \in \Lambda^M$:

$$(8.10) \quad \sum_{k=1}^{M} \theta_k \|\hat{\boldsymbol{\mu}}_k - g\|_2^2 = \|\hat{\boldsymbol{\mu}}_{\boldsymbol{\theta}} - g\|_2^2 + \text{pen}(\boldsymbol{\theta}).$$

This can be shown by using simple properties of the Euclidean norm, or by noting that the equality above is a bias-variance decomposition. For $g = 0$, (8.10) yields $\text{pen}(\boldsymbol{\theta}) = -\|\hat{\boldsymbol{\mu}}_{\boldsymbol{\theta}}\|_2^2 + \sum_{k=1}^{M} \theta_k \|\hat{\boldsymbol{\mu}}_k\|_2^2$.

LEMMA 8.2. *Let $F$ be any one of the functions $H_{\text{pen}}$, $V_{\text{pen}}$ or $W_{\text{pen}}$ defined in* (2.2), (5.1) *and* (6.6), *respectively. Then $F$ is convex, differentiable and satisfies for all $\boldsymbol{\theta}, \boldsymbol{\theta}_0 \in \Lambda^M$,*

$$(8.11) \quad F(\boldsymbol{\theta}) = F(\boldsymbol{\theta}_0) + \nabla F(\boldsymbol{\theta}_0)^T (\boldsymbol{\theta} - \boldsymbol{\theta}_0) + \frac{1}{2}\|\hat{\boldsymbol{\mu}}_{\boldsymbol{\theta}} - \hat{\boldsymbol{\mu}}_{\boldsymbol{\theta}_0}\|_2^2.$$

*Furthermore, if $\hat{\boldsymbol{\theta}}$ is a minimizer of $F$ over the simplex then for all $\boldsymbol{\theta} \in \Lambda^M$,*

$$(8.12) \quad F(\boldsymbol{\theta}) \ge F(\hat{\boldsymbol{\theta}}) + \frac{1}{2}\|\hat{\boldsymbol{\mu}}_{\boldsymbol{\theta}} - \hat{\boldsymbol{\mu}}_{\hat{\boldsymbol{\theta}}}\|_2^2.$$

PROOF. Using (8.10) with $g = 0$, we obtain that the function $F$ is a polynomial of degree 2, of the form $F(\boldsymbol{\theta}) = \text{affine}(\boldsymbol{\theta}) + \frac{1}{2}\|\hat{\boldsymbol{\mu}}_{\boldsymbol{\theta}}\|_2^2$ where $\text{affine}(\cdot)$ is an affine function of $\boldsymbol{\theta}$. This shows that $F$ is convex and differentiable. The result (8.11) follows by uniqueness of the Taylor expansion of $F$ [or by an explicit calculation of $\nabla F(\boldsymbol{\theta}_0)$]. Inequality (8.12) is a consequence of [8], 4.2.3, equation (4.21). □

8.5. *Lower bound.*

PROOF OF PROPOSITION 2.1. The lower bounds of [32], Theorem 5.4, are stated in expectation, but inspection of the proof of [32], Theorem 5.3 with $S = 1$, $\delta = \infty$ and $R = \log(1 + eM)$, reveals that the lower bound holds also in probability since it is an application of [38], Theorem 2.7. This result yields that there exist absolute constants $p, c, C > 0$ and $\mathbf{f}_1, \ldots, \mathbf{f}_M \in \mathbf{R}^n$ such that for any estimator $\hat{\boldsymbol{\mu}}$,

$$\sup_{j=1,\ldots,M} \mathbb{P}_{\mathbf{f}_j}(\Omega_j) \ge p, \qquad \Omega_j := \{\|\hat{\boldsymbol{\mu}} - \mathbf{f}_j\|_2^2 \ge c\sigma^2 \log(M)\},$$

provided that $\log(M) \le cn$ and $n, M > C$. Set $b_j = \mathbf{f}_j$ for all $j = 1, \ldots, M$. This lower bound implies that for any estimator $\hat{\boldsymbol{\mu}}$

$$\sup_{\mathbf{f} \in \mathbf{R}^n} \mathbb{P}_\mathbf{f}\left(\|\hat{\boldsymbol{\mu}} - \mathbf{f}\|_2^2 - \min_{k=1,\ldots,M} \|b_k - \mathbf{f}\|_2^2 \ge c\sigma^2 \log(M)\right) \ge p.$$

For all $j = 1, \ldots, M$, let $A_j = (1/\|\mathbf{f}_j\|_2^2)\mathbf{f}_j\mathbf{f}_j^T$ so that $A_j$ is the orthogonal projection on the linear span of $\mathbf{f}_j$. The orthoprojector $A_j$ has rank one so under $\mathbb{P}_{\mathbf{f}_j}$, $\|A_j\mathbf{y} - \mathbf{f}_j\|_2^2/\sigma^2$ is a $\chi^2$ random variable with one degree of freedom. Let $\Omega'_j$ be the event $\{\|A_j\mathbf{y} - \mathbf{f}_j\|_2^2 \le c\sigma^2 \log(M)/2\}$ and let $\bar{\Omega}'_j$ be its complementary event. A two-sided bound on the Gaussian tail implies that $\mathbb{P}_{\mathbf{f}_j}(\bar{\Omega}'_j) \le 2/(M^{c/4})$, which is smaller than $p/2$ if $M$ is larger than some absolute constant, so that we have $\mathbb{P}_{\mathbf{f}_j}(\bar{\Omega}_j \cup \bar{\Omega}'_j) \le 1 - p + p/2$ where $\bar{\Omega}_j$ is the complementary of $\Omega_j$, which implies $\mathbb{P}_{\mathbf{f}_j}(\Omega_j \cap \Omega'_j) \ge p/2$. Thus, for any estimator $\hat{\boldsymbol{\mu}}$ and $M$ large enough,

$$\sup_{\mathbf{f} \in \{\mathbf{f}_1,\ldots,\mathbf{f}_M\}} \mathbb{P}_\mathbf{f}\left(\|\hat{\boldsymbol{\mu}} - \mathbf{f}\|_2^2 - \min_{k=1,\ldots,M} \|A_k\mathbf{y} - \mathbf{f}\|_2^2 \ge c\sigma^2 \log(M)/2\right)$$
$$\ge p/2 =: p^*. \qquad \square$$

8.6. *Proof of Proposition* 3.1.

PROOF OF PROPOSITION 3.1. Let $a \in (0, 1)$. By definition of $\hat{J}$, we have for all $k = 1, \ldots, M$,

$$\|\hat{\boldsymbol{\mu}}_{\hat{J}} - \mathbf{f}\|_2^2 \le \|\hat{\boldsymbol{\mu}}_k - \mathbf{f}\|_2^2 + \Delta_{\hat{J}k} - \frac{a}{2}\|\hat{\boldsymbol{\mu}}_{\hat{J}} - \hat{\boldsymbol{\mu}}_k\|_2^2 + \frac{a}{2}\|\hat{\boldsymbol{\mu}}_{\hat{J}} - \hat{\boldsymbol{\mu}}_k\|_2^2$$
$$\le \|\hat{\boldsymbol{\mu}}_k - \mathbf{f}\|_2^2 + \frac{1}{a} \max_{j,k=1,\ldots,M}\left(a\Delta_{jk} - \frac{1}{2}\|a\hat{\boldsymbol{\mu}}_{\hat{J}} - a\hat{\boldsymbol{\mu}}_k\|_2^2\right)$$
$$+ a(\|\hat{\boldsymbol{\mu}}_{\hat{J}} - \mathbf{f}\|_2^2 + \|\mathbf{f} - \hat{\boldsymbol{\mu}}_k\|_2^2).$$

By rearranging, we have almost surely

$$(8.13) \qquad \|\hat{\boldsymbol{\mu}}_{\hat{J}} - \mathbf{f}\|_2^2 \le \frac{1+a}{1-a} \min_{k=1,\ldots,M} \|\hat{\boldsymbol{\mu}}_k - \mathbf{f}\|_2^2 + \frac{\Xi}{a(1-a)},$$

where

$$\Xi := \max_{j,k=1,\ldots,M}\left(2\boldsymbol{\xi}^T(\hat{\boldsymbol{\mu}}'_j - \hat{\boldsymbol{\mu}}'_k) - 2\sigma^2 \operatorname{Tr}(A'_j - A'_k) - \frac{1}{2}\|\hat{\boldsymbol{\mu}}'_j - \hat{\boldsymbol{\mu}}'_k\|_2^2\right),$$

and for all $j = 1, \ldots, M$, $\hat{\boldsymbol{\mu}}'_j := a\hat{\boldsymbol{\mu}}_j = A'_j\mathbf{y} + b'_j$, $A'_j := aA_j$, $b'_j := ab_j$. Simple algebra yields that

$$\frac{1}{2} \max_{j \ne k} \|A'_j - A'_k\|_2 \le a\phi \le \phi,$$

where $\phi$ is defined in (2.7). By Proposition 8.1, as in the proof of Theorem 2.1, we have $\Xi \leq 30\phi^2\sigma^2(x + 2\log M)$ with probability greater than $1 - \exp(-x)$.

Set $\varepsilon = 3a$ and choose the absolute constant $c > 0$ such that for all $\varepsilon < c$, $(1 + a)/(1 - a) \leq 1 + \varepsilon$ and $1/(1 - a) \leq 2$. □

8.7. *Smoothness adaptation.*

PROOF OF PROPOSITION 7.1. Because the ellipsoids are ordered, if $\mathbf{f} \in \Theta(\beta, L)$ then

$$\mathbb{E}\|\mathbf{f} - \tilde{\boldsymbol{\mu}}\|_2^2 \leq \min_{j:\beta_j \leq \beta} \mathbb{E}\|\mathbf{f} - A_{\beta_j}\mathbf{y}\|_2^2 + 60\sigma^2 \log M \leq \min_{j:\beta_j \leq \beta} C^* n^{\frac{1}{2\beta_j+1}}(1 + o(1)).$$

If $\beta \in [\beta_j, \beta_{j+1})$ for some $j$, then $\beta_{j+1} - \beta_j = \beta_j/(\log(n)\log\log n)$ and simple algebra yields

$$n^{\frac{1}{2\beta_j+1} - \frac{1}{2\beta+1}} \leq n^{\frac{2\beta_{j+1}-2\beta_j}{(2\beta+1)(2\beta_j+1)}} = n^{\frac{2\beta_j}{(2\beta+1)(2\beta_j+1)\log(n)\log\log n}}$$

$$\leq n^{\frac{1}{(2\beta+1)\log(n)\log\log n}} \leq e^{\frac{1}{3\log\log n}},$$

where we used that $\beta \geq 1$ for the last inequality.

Now assume that $\beta \geq \beta_M$. Let $\varepsilon_n = 1/(\log(n)\log\log n)$, and let $c = 120\log(1 + \varepsilon_3)/\varepsilon_3$. By definition of $M$,

$$\beta_M = e^{M\log(1 + \frac{1}{\log(n)\log\log n})} \geq e^{120\log\log(n)\frac{\log(1+\varepsilon_n)}{\varepsilon_n}}$$

$$\geq e^{c\log\log(n)} = \log(n)^c,$$

since the function $t \to \log(1+t)/t$ is decreasing and $n \geq 3$. A numerical approximation gives $c \geq 1.01$. Thus,

$$n^{\frac{1}{2\beta_M+1}}n^{\frac{-1}{2\beta+1}} \leq n^{\frac{1}{2\beta_M+1}} \leq e^{\frac{\log n}{2\beta_M}} \leq e^{\frac{1}{2\log(n)^{c-1}}}.$$

In summary, we have proved that $\min_{j:\beta_j \leq \beta} n^{\frac{1}{2\beta_j+1}} \leq n^{\frac{1}{2\beta+1}}(1 + o(1))$, thus

$$\sup_{\mathbf{f}\in\Theta(\beta,L)} \mathbb{E}\|\mathbf{f} - \tilde{\boldsymbol{\mu}}\|_2^2 \leq C^* n^{\frac{1}{2\beta+1}}(1 + o(1)). \quad \Box$$

8.8. *Convex aggregation.*

LEMMA 8.3 (Maurey argument). *Let $m$ and $\Lambda_m^M$ be defined in (7.2) and (7.3). Let $Q(\boldsymbol{\theta}) = \boldsymbol{\theta}^T \Sigma \boldsymbol{\theta} + \mathbf{v}^T \boldsymbol{\theta} + a$ for some semidefinite matrix $\Sigma$, $\mathbf{v} \in \mathbf{R}^M$ and $a \in \mathbf{R}$. Then*

(8.14) $$\min_{\boldsymbol{\theta}\in\Lambda_m^M} Q(\boldsymbol{\theta}) \leq \min_{\boldsymbol{\theta}\in\Lambda^M} Q(\boldsymbol{\theta}) + \frac{4\max_{j=1,\ldots,M}\Sigma_{jj}}{m}.$$

PROOF. Let $\boldsymbol{\theta}^* \in \Lambda^M \in \operatorname{argmin}_{\boldsymbol{\theta} \in \Lambda^M} Q(\boldsymbol{\theta})$. Let $\eta$ be a random variable valued in $\{\boldsymbol{e}_1, \ldots, \boldsymbol{e}_M\}$ such that $\mathbb{P}(\eta = \boldsymbol{e}_j) = \theta_j^*$ for all $j = 1, \ldots, M$, and let $\eta_1, \ldots, \eta_m$ be $m$ i.i.d. copies of $\eta$. The random variable $\bar{\eta} = \frac{1}{m} \sum_{q=1}^m$ is valued in $\Lambda_m^M$ and $\mathbb{E}\bar{\eta} = \boldsymbol{\theta}^*$. A bias variance decomposition and the independence of $\eta_1, \ldots, \eta_m$ yield

$$\mathbb{E} Q(\bar{\eta}) = Q(\boldsymbol{\theta}^*) + \frac{\mathbb{E}[(\eta_1 - \boldsymbol{\theta}^*)^T \Sigma (\eta_1 - \boldsymbol{\theta}^*)]}{m}.$$

Using the triangle inequality,

$$\mathbb{E}[(\eta_1 - \boldsymbol{\theta}^*)^T \Sigma (\eta_1 - \boldsymbol{\theta}^*)] \le 2(\boldsymbol{\theta}^*)^T \Sigma \boldsymbol{\theta}^* + 2\mathbb{E}[\eta_1^T \Sigma \eta_1] \le 4 \max_{j=1,\ldots,M} \Sigma_{jj}.$$

Since $\bar{\eta}$ is valued in $\Lambda_m^M$, $\min_{\boldsymbol{\theta} \in \Lambda_m^M} Q(\boldsymbol{\theta}) \le \mathbb{E} Q(\bar{\eta})$ and the proof is complete. □

PROOF OF (7.5) OF PROPOSITION 7.2. The condition on $M$, $n$ implies that $m \ge 1$ where $m$ is defined in (7.2). Let $C > 0$ be an absolute constant whose value may change from line to line. Applying Theorem 2.1 yields that on an event of probability greater than $1 - 2\exp(-x)$,

$$(8.15) \qquad \frac{1}{n} \|\hat{\boldsymbol{\mu}}_{\Lambda_m^M} - \mathbf{f}\|_2^2 \le \min_{\boldsymbol{\theta} \in \Lambda_m^M} \frac{1}{n} \|\hat{\boldsymbol{\mu}}_{\boldsymbol{\theta}} - \mathbf{f}\|_2^2 + \frac{C\sigma^2 (\log(|\Lambda_m^M|) + x)}{n}.$$

By [25], page 8, we have $\log |\Lambda_m^M| \le m \log \frac{2eM}{m}$. We use (8.14) with $Q(\boldsymbol{\theta}) = \|\hat{\boldsymbol{\mu}}_{\boldsymbol{\theta}} - \mathbf{f}\|_2^2$ to get

$$\min_{\boldsymbol{\theta} \in \Lambda_m^M} \frac{1}{n} \|\hat{\boldsymbol{\mu}}_{\boldsymbol{\theta}} - \mathbf{f}\|_2^2 \le \min_{\boldsymbol{\theta} \in \Lambda^M} \frac{1}{n} \|\hat{\boldsymbol{\mu}}_{\boldsymbol{\theta}} - \mathbf{f}\|_2^2 + \frac{4}{nm} \max_{j=1,\ldots,M} \|\hat{\boldsymbol{\mu}}_j\|_2^2.$$

We have $(1/n) \max_{j=1,\ldots,M} \|\hat{\boldsymbol{\mu}}_j\|_2^2 \le C(\|\boldsymbol{\xi}\|_2^2/n + R^2) \le C(\sigma^2(2 + 3x) + R^2)$ on an event of probability at least $1 - \exp(-x)$, where for the second inequality we used (3.5) with $B = I_{n \times n}$. Thus, with probability greater than $1 - e^{-x}$,

$$(8.16) \qquad \min_{\boldsymbol{\theta} \in \Lambda_m^M} \frac{1}{n} \|\hat{\boldsymbol{\mu}}_{\boldsymbol{\theta}} - \mathbf{f}\|_2^2 \le \min_{\boldsymbol{\theta} \in \Lambda^M} \frac{1}{n} \|\hat{\boldsymbol{\mu}}_{\boldsymbol{\theta}} - \mathbf{f}\|_2^2 + \frac{C(\sigma^2(2 + 3x) + R^2)}{m}.$$

Simple algebra yields that

$$(8.17) \qquad \begin{aligned} \frac{1}{m} &\le C\sqrt{\frac{\log(1 + M/\sqrt{n})}{n}}, \\ \frac{m \log(2eM/m)}{n} &\le C\sqrt{\frac{\log(1 + M/\sqrt{n})}{n}}. \end{aligned}$$

Combining (8.15), (8.16) and (8.17) with the union bound completes the proof. □

PROOF OF (7.4) OF PROPOSITION 7.2. Let $\boldsymbol{\theta} \in \Lambda^M$. By definition of $\hat{\boldsymbol{\theta}}_C$, $C_p(\hat{\boldsymbol{\theta}}) \le C_p(\boldsymbol{\theta})$. This can be rewritten as

$$\|\hat{\boldsymbol{\mu}}_{\hat{\boldsymbol{\theta}}_C} - \mathbf{f}\|_2^2 \le \|\hat{\boldsymbol{\mu}}_{\hat{\boldsymbol{\theta}}} - \mathbf{f}\|_2^2 + 2\boldsymbol{\xi}^T(\hat{\boldsymbol{\mu}}_{\hat{\boldsymbol{\theta}}} - \hat{\boldsymbol{\mu}}_{\boldsymbol{\theta}}) - 2\sigma^2 \operatorname{Tr}(A_{\hat{\boldsymbol{\theta}}} - A_{\boldsymbol{\theta}}).$$

The function $(\boldsymbol{\theta}', \boldsymbol{\theta}) \to 2\boldsymbol{\xi}^T(\hat{\boldsymbol{\mu}}_{\hat{\boldsymbol{\theta}}} - \hat{\boldsymbol{\mu}}_{\boldsymbol{\theta}}) - 2\sigma^2 \operatorname{Tr}(A_{\hat{\boldsymbol{\theta}}} - A_{\boldsymbol{\theta}})$ is linear in $\boldsymbol{\theta}$ and linear in $\hat{\boldsymbol{\theta}}$, thus it is maximized at vertices of the simplex and

$$2\boldsymbol{\xi}^T(\hat{\boldsymbol{\mu}}_{\hat{\boldsymbol{\theta}}} - \hat{\boldsymbol{\mu}}_{\boldsymbol{\theta}}) - 2\sigma^2 \operatorname{Tr}(A_{\hat{\boldsymbol{\theta}}} - A_{\boldsymbol{\theta}}) \le \max_{j,k=1,\ldots,M} 2\boldsymbol{\xi}^T(\hat{\boldsymbol{\mu}}_k - \hat{\boldsymbol{\mu}}_j) - 2\operatorname{Tr}(A_j - A_k)$$

$$= \max_{j,k=1,\ldots,M} \Delta_{jk},$$

where $\Delta_{jk}$ is defined in (3.3). Fix some pair $(j, k)$. Let $B = A_j - A_k$ and $\boldsymbol{b} = (A_j - A_k)\mathbf{f} + \boldsymbol{b}_j - \boldsymbol{b}_k$. We have $\|\|B\|\|_2 \le 2$, $\|B\|_F \le \|\|B\|\|_2 \|I_{n \times n}\|_F \le 2\sqrt{n}$ and $\|\boldsymbol{b}\|_2 \le 4R\sqrt{n}$. We apply (3.5) to the matrix $B$ and (3.4) to the vector $\boldsymbol{b}$, which yields that with probability greater than $1 - 2\exp(-x)$,

$$\Delta_{jk} \le 8(\sigma^2 + \sigma R\sqrt{2})\sqrt{nx} + 8\sigma^2 x.$$

The union bound over all pairs $j, k = 1, \ldots, M$ completes the proof. □

**Acknowledgements.** The author would like to thank the Associate Editor for valuable comments that improved the paper substantially. The author also thanks Alexandre Tsybakov for some helpful discussions on early versions of this manuscript.

DEPARTMENT OF STATISTICS & BIOSTATISTICS
RUTGERS, THE STATE UNIVERSITY OF NEW JERSEY
501 HILL CENTER, BUSCH CAMPUS
110 FRELINGHUYSEN ROAD
PISCATAWAY, NEW JERSEY 08854
USA
E-MAIL: pcb71@stat.rutgers.edu